\documentclass[a4paper,11pt]{amsart}

\usepackage{amssymb}
\usepackage{amscd}

\newcommand{\nn}[1]{(\ref{#1})}

\newcommand{\tth}{\h}

\newcommand{\tM}{{\tilde M}}
\newcommand{\tQ}{{\tilde Q}}
\newcommand{\tH}{{\tilde H}}
\newcommand{\tW}{{\tilde W}}

\def\frak{\mathfrak}
\def\Bbb{\mathbb}
\def\Cal{\mathcal}

\newtheorem*{prop*}{Proposition
\thesubsection}
\newtheorem*{thm*}{Theorem \thesubsection}
 \newtheorem*{lem*}{Lemma
\thesubsection} 
\newtheorem*{cor*}{Corollary \thesubsection}

\newcommand{\ol}[1]{\overline{#1}}
\newcommand{\wh}[1]{\widehat{#1}}

\newcommand{\alb}{{\overline{\alpha}}}
\newcommand{\beb}{{\overline{\beta}}}

\renewcommand{\k}{{\mathbf{k}}}

\newcommand{\id}{\operatorname{id}}

\let\ker\Ker

\newcommand{\Le}{{\Cal L^{\theta}}}

\newcommand{\x}{\times}
\renewcommand{\o}{\circ}

\def\crn#1#2{{\vcenter{\vbox{
        \hbox{\kern#2pt \vrule width.#2pt height#1pt
           }
          \hrule height.#2pt}}}}

\let\ccdot\cdot
\def\cdot{\hbox to 2.5pt{\hss$\ccdot$\hss}}

\newcommand{\al}{\alpha}
\newcommand{\be}{\beta}
\newcommand{\ga}{\gamma}

\newcommand{\ka}{\kappa}
\newcommand{\la}{\lambda}
\newcommand{\om}{\omega}
\renewcommand{\phi}{\varphi}
\newcommand{\ph}{\varphi}

\newcommand{\si}{\sigma}
\renewcommand{\th}{\theta}
\newcommand{\ze}{\zeta}

\newcommand{\La}{\Lambda}
\newcommand{\Ga}{\Gamma}
\newcommand{\Ph}{\Phi}
\newcommand{\Om}{\Omega}

\newcommand{\lpl}                         
{\mbox{$
\begin{picture}(12.7,8)(-.5,-1)
\put(2,0.2){$+$}
\put(6.2,2.8){\oval(8,8)[l]}
\end{picture}$}}

\newcommand{\fg}{{\frak g}}

\newcommand{\tg}{{\tilde {\frak g}}}

\newcommand{\tp}{{\tilde {\frak p}}}

\newcommand{\ce}{{\Cal E}}
\newcommand{\ct}{{\Cal T}}
\newcommand{\tca}{{\tilde{\Cal A}}}
\newcommand{\tce}{{\tilde{\Cal E}}}
\newcommand{\tcg}{{\tilde{\Cal G}}}
\newcommand{\tcs}{{\tilde{\Cal S}}}
\newcommand{\tct}{{\tilde{\Cal T}}}
\newcommand{\tcw}{{\tilde{\Cal W}}}

\newcommand{\bZ}{{\Bbb Z}}

\newcommand{\bg}{{\boldsymbol g}}
\newcommand{\bh}{{\boldsymbol h}}

\newcommand{\nd}{\nabla}
\newcommand{\Ps}{\Psi}
\newcommand{\Rho}{P}
\newcommand{\Up}{\Upsilon}

\newcommand{\K}{{\sf K}}
\renewcommand{\L}{{\sf L}}

\def\sideremark#1{\ifvmode\leavevmode\fi\vadjust{\vbox to0pt{\vss
 \hbox to 0pt{\hskip\hsize\hskip1em
 \vbox{\hsize3cm\tiny\raggedright\pretolerance10000
 \noindent #1\hfill}\hss}\vbox to8pt{\vfil}\vss}}}%
\begin{document}
\title[CR--tractors and the Fefferman space]{CR--Tractors and the
Fefferman Space}  
\author{Andreas \v Cap and A.\ Rod Gover}
\date{November 30, 2006}
\address{A.C.: Institut f\"ur Mathematik, Universit\"at Wien, Strudlhofgasse 4,
A--1090 Wien, Austria and International Erwin Schr\"odinger Institute for
Mathematical Physics, Boltzmanngasse 9, A--1090 Wien, Austria\newline\indent
A.R.G.: Department of Mathematics, The University of Auckland, Auckland, 
New Zealand}
\email{Andreas.Cap@esi.ac.at, gover@math.auckland.ac.nz}
\subjclass{primary:32V05, 53B15, 53C15, 53C25, 58J60, 58J70; secondary:32V20, 53B30, 53C29}
\keywords{CR structure, Fefferman space, conformal structure, tractor calculus,
invariant differential operator}
\begin{abstract} 
  We develop the natural tractor calculi associated to conformal and
  CR structures as a fundamental tool for the study of Fefferman's
  construction of a canonical conformal class on the total space of a
  circle bundle over a non--degenerate CR manifold of hypersurface
  type. In particular we construct and treat the basic objects that
  relate the natural bundles and natural operators on the two spaces.
  This is illustrated with several applications: We prove that a
  number of conformally invariant overdetermined systems, including
  Killing form equations and the equations for twistor spinors, admit
  non--trivial solutions on any Fefferman space. We show that, on a
  Fefferman space, the space of infinitesimal conformal isometries
  naturally decomposes into a direct sum of subspaces, which admit an
  interpretation as solutions of certain CR invariant PDE's. Finally
  we explicitly analyse the relation between tractor calculus on a CR
  manifold and the complexified conformal tractor calculus on its
  Fefferman space, thus obtaining a powerful computational tool for
  working with the Fefferman construction.
\end{abstract}

\maketitle

\subsection*{Note}: This article has been published in Indiana Univ. 
Math. J. \textbf{57} No. 5 (2008) 2519-2570 with (apart from other small changes) 
a different numbering of statements. The published version is acessible online 
via \newline
http://www.iumj.indiana.edu/IUMJ/fulltext.php?year=2008\&volume=57\&artid=3359 

\section{Introduction}\label{1}
Given a non--degenerate CR manifold $M$ of hypersurface type, the
Fefferman space $\tilde M$, associated to $M$, is the total space of a
circle bundle over $M$, on which the given CR structure induces a
natural (indefinite) conformal structure. Fefferman's original construction in
\cite{Fefferman} applied to smooth boundaries of strictly pseudoconvex
domains and it was a vision of Fefferman to exploit this relationship
to study the invariant theory and geometry of these objects in terms
of conformal geometry. Despite the subsequent and increasing interest
in these, and related issues of CR geometry, this remarkable
construction has not been fully explored.  The explanation surely lies
in the complicated relationship between natural bundles and operators
on the two structures. The main aim of this work is to present a
treatment of the Fefferman space that provides both a conceptual and
calculationally practical solution to this problem at the level of the
basic foundational issues. We illustrate the power of the approach by
recovering very quickly a string of applications that we discuss
below.

The Fefferman construction was generalised to abstract CR manifolds in
\cite{BDS} using the canonical Chern--Moser Cartan connection. A main
application of this and \cite{Fefferman} was that chains in $M$ can be
interpreted as the projections of null geodesics in $\tilde M$. This
crucially depends on the fact that the canonical Cartan connection of
the conformal structure on $\tilde M$ is closely related to the CR
Cartan connection on $M$. In both \cite{BDS} and \cite{Fefferman} the
necessary verification of this relationship was evidently very
complicated, with many details of the calculation sketched or omitted.

A different approach to the Fefferman construction was developed by
J.~Lee in \cite{LeeF}. For a choice of a pseudo--Hermitian structure
(i.e.~a contact form) on $M$, Lee used the associated Webster--Tanaka
connection, to directly define a metric on the Fefferman space $\tilde
M$. Then he showed that changing the contact form only leads to a
conformal rescaling of that metric. A feature of this approach is that
it directly leads to explicit formulae. On the other hand, the
relation between the Cartan connections is not directly visibly in
this picture.  Using Lee's approach, it was realised that Fefferman
spaces provide particularly interesting examples of conformal
structures. For example, it was shown in \cite{Lewandowski} for
CR--dimension one, and in \cite{Baum} in general, that (slightly
modified) Fefferman spaces always admit non--trivial twistor spinors.

Conformal and CR geometry have been very active areas of research
recently. Questions related to invariant operators with symbol a power
of the Laplacian, respectively the sub-Laplacian, as well as to
Branson's $Q$--curvature have received particular interest. These
studies are driven not only by differential geometry but also through
the role of these objects in geometric analysis. It is not surprising
then that there is renewed interest in Fefferman's construction as a
bridge between CR geometry and conformal geometry. In \cite{FeffHir},
via the Fefferman metric, a CR Q-curvature was defined and studied,
while in \cite{Gover-Graham} the Fefferman metric was used as one of
the two main construction techniques for the CR invariant powers of
the sub-Laplacian.  In these treatments the use of the Fefferman
structure is relatively straightforward since the central objects
involved are density bundles and operators between them. To develop
analogous results involving, for example, tensor or spinor bundles
requires a significantly deeper understanding of the Fefferman space
and its precise geometric relationship to the underlying CR structure.
The problem is that the relation between such bundles on a Fefferman
space and irreducible bundles on the underlying CR manifold is very
complicated in general. 

We solve this problem by introducing tractor bundles and tractor
connections as a new tool in the study of the Fefferman construction; 
tensor and spinor bundles are captured as subquotients in these. 
The tractor bundles are natural vector bundles associated to a CR strucure
or a conformal structure (and more generally to a so--called parabolic
geometry) which are endowed with canonical invariant linear
connections.  They are equivalent to the Cartan princpal bundle and
the canonical Cartan connection, see \cite{TAMS}, but are easier to
handle. We show that the relation between conformal tractor bundles on
a Fefferman space $\tilde M$ and CR tractor bundles on the underlying
CR manifold $M$ is rather simple. For example, CR standard tractors on
$M$ may be identified with conformal standard tractors on $\tilde M$ which are
parallel for the canonical tractor connection in the direction of the
fibers of $\tilde M\to M$. This is  simpler than the relation
between the canonical Cartan geometries, which involves extension of
the structure group and equivariant extension of the Cartan
connection. Irreducbile bundles then naturally occur as subquotients
of tractor bundles, so we obtain a vehicle for effectively carrying
the conceputal and calculational details of the relationship.

To have tractors at our disposal, in addition to the CR structure on
$M$, we have to choose a certain root of the canonical bundle, compare
with \cite{Gover-Graham}. This is no restriction locally for arbitrary
structures, and is no restriction globally in the embedded case. This
slightly richer structure leads to an immediate payoff: we
automatically get a canonical spin structure on the Fefferman space,
and we exactly recover Fefferman's original construction (which is an
$(n+2)$--fold covering of the one in \cite{BDS}) for embedded
structures. Moreover, using Tanaka's version of the canonical Cartan
connection, the construction automatically extends to the class of
partially integrable almost CR structures, which is much larger than
the integrable ones. While a canonical conformal structure is obtained
without the integrability assumption, the close relation between the
canonical Cartan connections surprisingly is true only in the
integrable case. We prove this in the language of tractor connections
in Theorem \ref{2.5}.

Having this at hand, we apply the powerful tools available for
parabolic geometries to study the relation between a CR manifold and
its Fefferman space as well as the conformal geometry of Fefferman
spaces in section \ref{3}. We obtain a short and conceptual proof of
the existence of non--trivial twistor spinors, new results on the
existence and construction of odd degree conformal Killing forms, as
well as a natural decomposition of conformal Killing vector fields
(i.e.\ infinitesimal conformal isometries). A crucial ingredient in
all this is that on the conformal standard tractor bundle of a
Fefferman space one obtains a parallel, orthogonal complex
structure. In a companion article \cite{characterization} to this one
we show that Fefferman spaces are characterised (up to local conformal
isometry) by the existence of such a complex structure, which can be
viewed as a restriction on the so--called conformal holonomy group.
This also leads to a new proof and extension of Sparling's
characterisation of Fefferman spaces from \cite{Graham}.

In section \ref{4}, the abstract developments from preceeding sections
are converted into an explicit calculus. We first show that the
tractor calculus developed in \cite{Gover-Graham} recovers precisely
the complexified normal CR standard tractor bundle. This result should
be of independent interest. Next, we study the consequences of the
existence of a parallel, orthogonal complex structure in terms of
conformal tractor calculus. In the case of a Fefferman space $\tilde
M\to M$, we obtain explicit relations between CR tractor calculus on
$M$ and conformal tractor calculus on $\tilde M$. Using this, we
completely decsribe how a contact/pseudo--Hermitian structure on the
CR manifold may be related to an equivalent structure on the Fefferman
space. We illustrate the utility of this by converting the
developments of section \ref{3} into explicit formulae. Moreover, we
explicitly compute the metric in the conformal class associated to a
choice of CR scale, thus tying in with Lee's approach. We also obtain
explicit relations between Webster--Tanaka connections on $M$ and
Levi--Civita connections on $\tilde M$. As an application, we discuss
a tractor interpretation of Einstein-type structures (cf.\
\cite{Leitner:TSPE}) in CR geometry.

\subsection*{Acknowledgement} First author supported by projects
P15747-N05 and P195900-N13 of the Fonds zur F\"orderung der
wissenschaftlichen Forschung (FWF). The second author would like to
thank the Royal Society of New Zealand for support via Marsden Grant
no.\ 02-UOA-108, and the New Zealand Institute of Mathematics and its
Applications for support via a Maclaurin Fellowship.

\section{The Fefferman space}\label{2}

\subsection{CR manifolds}\label{2.1}
An \textit{almost CR structure} of hypersurface type on a smooth
manifold $M$ of dimension $2n+1$ is a rank $n$ complex subbundle $H$
of the tangent bundle $TM$. We denote by $J:H\to H$ the complex
structure on the subbundle. The quotient $Q:=TM/H$ is a real line
bundle on $M$. Let $q:TM\to Q$ be the obvious surjection. For two
sections $\xi,\eta\in\Ga(H)$ the expression $q([\xi,\eta])$ is
bilinear over smooth functions, and so we obtain a skew symmetric bundle
map $\Cal L:H\x H\to Q$ given by $\Cal
L(\xi(x),\eta(x))= q([\xi,\eta](x))$. The almost CR structure is called
\textit{non--degenerate} if $\Cal L(\xi,\eta)=0$ for all $\eta$
implies $\xi=0$. This is equivalent to the fact that $H$, viewed as a
real subbundle of $TM$, defines a contact structure on $M$.

Looking at the complexified tangent bundle $T_{\Bbb C}M$, the complex
structure on $H$ is equivalent to a splitting of the subbundle
$H_{\Bbb C}$ into the direct sum of the holomorphic part $H^{1,0}$ and
the antiholomorphic part $H^{0,1}=\overline{H^{1,0}}$. The almost CR
structure is called \textit{integrable} or a \textit{CR structure} if
the subbundle $H^{1,0}\subset T_{\Bbb C}M$ is involutive, i.e.~the
space of its sections is closed under the Lie bracket. A weakening of
this condition, called \textit{partial integrability}, is obtained by
requiring that the bracket of two sections of $H^{1,0}$ is a section
of $\Bbb H_{\Bbb C}$. This is equivalent to to $\Cal L$ being of type
$(1,1)$, that is $\Cal L(J\xi,J\eta)=\Cal L(\xi,\eta)$ for all
$\xi,\eta\in H$. Throughout the paper, we will only deal with
non--degenerate partially integrable CR structures, and unless
explicitly specified we will assume integrability.

Let $Q_{\Bbb C}$ be the complexification of $Q$ and let $q_{\Bbb
  C}:T_{\Bbb C}M\to Q_{\Bbb C}$ be the complex linear extension of $q$. The
\textit{Levi form} $\Cal L_{\Bbb C}$ of an almost CR structure is the
$Q_{\Bbb C}$--valued Hermitian form on $H^{1,0}$ induced by
$(\xi,\eta)\mapsto 2iq_{\Bbb C}([\xi,\overline{\eta}])$. Assuming
partial integrability, $\Cal L$ can be naturally identified with the
imaginary part of $\Cal L_{\Bbb C}$, and so under this assumption,
non--degeneracy of the almost CR structure also can be characterised
by non--degeneracy of the Levi form.

Choosing a local trivialisation of $Q$ and using the induced
trivialisation of $Q_{\Bbb C}$, $\Cal L_{\Bbb C}$ gives rise to a
Hermitian form. If $(p,q)$ is the signature of this form, then one
also says that $M$ is non--degenerate of signature $(p,q)$. If $p\neq
q$, then such local trivialisations of $Q$ necessarily fit
together to give a global trivialisation. In the case of symmetric
signature $(p,p)$ we assume that a global trivialisation of $Q$
exists. A global trivialisation of $Q$ is equivalent to a ray subbundle
of the line bundle of contact forms for $H\subset TM$, so we obtain a
notion of positivity for contact forms. 

Generic real hypersurfaces in complex manifolds give the prototypal
examples of CR structures, and form an important class for many
applications. Let $\Cal M$ be a complex manifold of complex dimension
$n+1$ and let $M\subset\Cal M$ be a smooth real hypersurface. For each
point $x\in M$ the tangent space $T_xM$ is a subspace of the complex
vector space $T_x\Cal M$ of real codimension one. This implies that
the maximal complex subspace $H_x$ of $T_xM$ must be of complex
dimension $n$. Of course, these subspaces fit together to define a
smooth subbundle $H\subset TM$, which by construction is equipped with
a complex structure. Since the bundle $H^{1,0}\subset T_{\Bbb C}M$ can
be viewed as the intersection of the involutive subbundles $T_{\Bbb
  C}M$ and $T^{1,0}\Cal M$ of $T_{\Bbb C}\Cal M|_M$ we see that we
always obtain a CR structure in this way. This structure is
non--degenerate if $H$ defines a contact structure on $M$, which is
satisfied generically. Examples of this type are usually referred to
as \textit{embedded CR manifolds}, in particular for $\Cal M=\Bbb
C^{n+1}$.

\subsection{The Fefferman space for the homogeneous model}\label{2.2}
For the case of the homogeneous model, the Fefferman construction can
be easily described and this motivates the general construction. 
Fix a complex vector space $\Bbb V$ of dimension
$n+2$, endowed with a Hermitian inner product $\langle\ ,\ \rangle$ of
signature $(p+1,q+1)$, where $p+q=n$. Let $\Cal C\subset\Bbb V$ be the
cone of nonzero null vectors, and let $M$ be the image of $\Cal C$ in
the complex projectivisation $\Cal P\Bbb V\cong\Bbb CP^{n+1}$. Hence
$M$ is the space of those complex lines in $\Bbb V$ which are null
with respect to $\langle\ ,\ \rangle$. Note that $M$ is a
smooth real hypersurface in $\Cal P\Bbb V$. The resulting CR structure
on $M$ is easily described explicitly: Given a null line
$\ell\subset\Bbb V$, the CR subspace $H_\ell M\subset T_\ell M$ is the
image of the complex orthogonal complement $\ell^\perp$ of $\ell$
under the tangent map of the obvious projection $\Cal C\to M$. Hence
we obtain an isomorphism $\ell^\perp/\ell\to H_\ell M$. It is easy to
see that under this isomorphism the Levi form, at $\ell$, corresponds (up to a
nonzero multiple) to the Hermitian form on $\ell^\perp/\ell$ induced by
$\langle\ ,\ \rangle$ in the obvious way. In particular, $M$ is
non--degenerate of signature $(p,q)$.

Let $G$ be the special unitary group of $(\Bbb V,\langle\ ,\ 
\rangle)$.  The standard linear action of $G$ on $\Bbb V$ restricts to
an action on $\Cal C$ and then descends to a smooth left action of $G$
on $M$.  Since the CR structure on $M$ is completely determined by the
Hermitian form $\langle\ ,\ \rangle$ it is evident the $G$ acts on $M$
by CR automorphisms. Elementary linear algebra shows that $G$ acts
transitively on $\Cal C$ and thus also on $M$.  Fixing a null line
$\ell\subset V$ and denoting by $P\subset G$ the stabiliser of $\ell$
we obtain an identification of $M$ with the homogeneous space $G/P$.

Note that the action of $G$ on $M$ is not effective. The kernel of the
action is the centre $Z(G)$ which consists of those multiples of the
identity which lie in $G$. Hence the possible factors are the
$(n+2)$nd roots of unity, and $Z(G)\cong\Bbb Z_{n+2}$. It is a
classical result that the action on $M$ induces an isomorphism between
$G/Z(G)$ and the group of CR automorphisms of $M$. So
$M=(G/Z(G))/(P/Z(G))$ as a homogeneous space for its group of CR
automorphisms.

To keep track of the full group $G$ we incorporate the principal
$\Bbb C^*$--bundle $\Cal C\to M$ as a part of the structure, since
the lift of the $G$--action to this bundle separates the points in the
centre. We shall work with this richer structure.

The Fefferman space $\tilde M$ arises naturally from the
underlying real picture. Let $\Cal P_{\Bbb R}\Bbb V$ be the real
projectivisation of $\Bbb V$ and let $\tilde M$ be the image of $\Cal
C$ in this real projective space. That is $\tilde M$ is the space of all
real lines in $\Bbb V$ which are null for the inner product $\langle\
,\ \rangle_{\Bbb R}$, the real part of $\langle\ ,\ \rangle$. The
space $\tilde M$ is a smooth hypersurface in $\Bbb RP^{2n+3}$, and we
have an obvious projection $\Cal C\to\tilde M$, which is a principal
bundle with fibre group $\Bbb R^*$. 

Any real null line generates a complex null line containing it.  This
smooth projection $\tilde M\to M$ is a 
fibre bundle over $M$,
with fibre the space $\Bbb RP^1\cong S^1$ of real lines in $\Bbb C$.

Fixing an element $v\in\Cal C$ and denoting by $\tilde\ell$ the real
line $\Bbb Rv$, the tangent map in $v$ of the projection $\Cal
C\to\tilde M$ identifies $T_{\tilde\ell}\tilde M$ with the quotient of
the real orthocomplement of $\tilde{\ell}$ by $\tilde\ell$. The inner
product $\langle\ ,\ \rangle_{\Bbb R}$ thus induces an inner product
on $T_\ell\tilde M$ and changing the point $v$ leads to a positive
rescaling of this product. Hence we obtain a well--defined conformal
structure on $\tilde M$.

Let $\tilde G$ be the connected component of the identity of the
orthogonal group of $(\Bbb V,\langle\ ,\ \rangle_{\Bbb R})$, and let
$\tilde P\subset\tilde G$ be the stabiliser of a real null line. Then
as above we obtain a transitive action of $\tilde G$ on $\tilde M$
which leads to an identification $\tilde M\cong\tilde G/\tilde P$. By
construction, $\tilde G$ acts by conformal isometries on $\tilde M$.
It is a classical result, this action identifies $\tilde G/Z(\tilde
G)$ with the group of conformal isometries of $\tilde M$.

We have noted above that the subgroup $G\subset\tilde G$ acts
transitively on $\Cal C$, so it also acts transitively on $\tilde M$.
Taking a real null line and the complex null line generated by it as
the base points of $\tilde M$ and $M$, we see that $G\cap\tilde
P\subset P$, and $G\cap\tilde P$ is the stabiliser of a real null
line, so $\tilde M\cong G/(G\cap\tilde P)$.

The inclusion $G\hookrightarrow\tilde G$ may be viewed as inducing a
reduction of structure group of the bundle $\tilde G\to\tilde M$ from
$\tilde P$ to $G\cap\tilde P\subset P$. This reduction is determined
by the complex structure on $\Bbb V$, which is orthogonal for the
inner product $\langle\ ,\ \rangle_{\Bbb R}$. Equivalently, we can
view this complex structure as a splitting $\Bbb V\otimes\Bbb C=\Bbb
V^{1,0}\oplus\Bbb V^{0,1}$, of the complexification of $\Bbb V$, into
a holomorphic and an anti--holomorphic part. The inner product
$\langle\ ,\ \rangle_{\Bbb R}$ induces a non-degenerate complex
bilinear form on $\Bbb V\otimes\Bbb C$. Since the complex structure is
orthogonal for $\langle\ ,\ \rangle_{\Bbb R}$, the subspaces $\Bbb
V^{1,0}$ and $\Bbb V^{0,1}$ are isotropic for this complex inner
product, while at the same time the induced pairing between the two
spaces exactly comes from viewing $\langle\ ,\ \rangle$ as a complex
bilinear map $\Bbb V\x\overline{\Bbb V}\to\Bbb C$.

\subsection{The canonical Cartan connection}\label{2.3}
For an arbitrary non--degenerate CR manifold $(M,H)$, the canonical
Cartan connection gives a description which is appropriate for
generalising the above construction for the homogeneous space
$(G/Z(G))/(P/Z(G))$. Various constructions of the canonical
Cartan bundle and Cartan connection can be found in
\cite{Chern-Moser,Tan76,Morimoto,Kuranishi1,CS}. The outcome may be
described as follows: First, one builds a principal fibre
bundle $p:\underline{\Cal G}\to M$ with structure group $P/Z(G)$;
this may be obtained as an extension of an adapted frame bundle of
$H\to M$ or as a subbundle of the (co)frame bundle of the total space
of the line bundle $Q\to M$. The principal bundle can be endowed with
a \textit{Cartan connection} $\underline{\om}\in\Om^1(\underline{\Cal
  G},\frak g)$, where $\frak g$ denotes the Lie algebra of $G$. This
generalises the trivialisation of the tangent bundle of $G/Z(G)$ by
left translations. 

Explicitly, we require that $\underline{\om}$ defines a trivialisation
of $T\underline{\Cal G}$, which is $P/Z(G)$--equivariant and
reproduces the generators of fundamental vector fields. If one
requires the curvature of $\underline{\om}$ to satisfy a normalisation
condition, which will be discussed in detail below, then the pair
$(\underline{\Cal G},\underline{\om})$ is uniquely determined up to
isomorphism.

For later purposes, it will be very important to extend the structure
group from $P/Z(G)$ to $P$. In the case of the homogeneous model, the
way to expose the centre is via its action on the restriction of the
tautological bundle to the hyperquadric. To generalise this to
arbitrary CR manifolds, one proceeds as follows: The natural choice of
a complex line bundle on a CR manifold is provided by the canonical
bundle $\Cal K$. By definition, $\Cal K$ is the $(n+1)$st complex
exterior power of the annihilator of $H^{0,1}$ in the complexified
cotangent bundle. In the case of the homogeneous model one shows that
$\Cal K$ is associated to the $\Bbb C^*$--bundle $\Cal C\to G/P$ with
respect to the representation $z\mapsto z^{-n-2}$, so the null cone
may be naturally viewed as the dual of a $(n+2)$nd root of the
canonical bundle.

When dealing with a general CR manifold $M$, we will always assume
that we have chosen a complex line bundle $\Cal E(1,0)\to M$ together
with a duality between $\Cal E(1,0)^{\otimes^{(n+2)}}$ and the
canonical bundle $\Cal K$. In general, such a choice may not exist
globally but locally it poses no problem. Moreover, for CR manifolds
embedded in $\Bbb C^{n+1}$ the canonical bundle is trivial, so the
required identification exists globally in this setting. For
$w,w'\in\Bbb R$ such that $w'-w\in\Bbb Z$, the map $\la\mapsto
|\la|^{2w}\overline{\la}^{(w'-w)}$ is a well defined one--dimensional
representation of $\Bbb C^*$. Hence we can define a complex line
bundle $\Cal E(w,w')$ over $M$ by forming the associated bundle to the
frame bundle of $\Cal E(1,0)$ with respect to this representation. By
construction we get $\Cal E(w',w)=\overline{\Cal E(w,w')}$, $\Cal
E(-w,-w')=\Cal E(w,w')^*$ and $\Cal E(k,0)=\Cal E(1,0)^{\otimes^k}$
for $k\in\Bbb N$. Finally, by definition $\Cal K\cong\Cal E(0,-n-2)$
in this notation. 

There is a natural inclusion of the real line bundle $Q:=TM/H$ into
the density bundle $\ce(1,1)$ which is defined as follows. For a local
nonzero section $\al$ of $\ce(1,0)$ one can, by definition, view
$\al^{-(n+2)}$ as a section of the canonical bundle $\Cal K$. Then by
\cite[Lemma 3.2]{LeeF}, there is a unique positive contact form $\th$,
with respect to which $\al^{-(n+2)}$ is length normalised. Mapping
$\xi\in TM$ to $\th(\xi)\al\overline{\al}$ then descends to an
inclusion, which is CR invariant.

Recall that the Cartan connection on $\underline{\Cal G}\to M$
identifies $TM$ with a bundle associated to the Cartan bundle
$\underline{\Cal G}\to M$, see e.g.~\cite[2.8]{TAMS}. Thus, also the
dual $\Cal K^*$ of $\Cal K$ is associated to $\underline{\Cal G}$, so
its frame bundle is a quotient of $\underline{\Cal G}$. By the
definition of $\Cal E(1,0)$, its frame bundle is an $(n+2)$--fold
covering of the frame bundle of $\Cal K^*$. Pulling back the covering
to $\underline{\Cal G}$, we obtain an $(n+2)$--fold covering $\Cal G$
of $\underline{\Cal G}$. The principal action of $P/Z(G)$ on
$\underline{\Cal G}$ lifts to a right action of $P$ on $\Cal G$. This
can be used to make the bundle $\Cal G\to M$ into a principal
$P$--bundle in such a way that $\Cal E(1,0)=G\x_P(\Bbb V^1)^*$, where
$\Bbb V^1\subset\Bbb V$ is the null line stabilised by $P$. A normal
Cartan connection $\om\in\Om^1(\Cal G,\frak g)$ on $\Cal G$ is
obtained by pulling back $\underline{\om}$.  As with the structure
$(\underline{\Cal G},\underline{\om})$, the pair $(\Cal G,\om)$ is
determined by $(M,H,\Cal E(1,0))$ up to isomorphism. We will show in
\ref{4.3} below that the Cartan structure $(\Cal G,\om)$ may be
recovered from the tractor bundle introduced in \cite{Gover-Graham}.

\subsection{The Fefferman space}\label{2.4}
Suppose we have given a partially integrable almost CR manifold
$(M,H)$ with a fixed choice $\ce(1,0)$ of an $(n+2)$nd root of the
anticanonical bundle as in \ref{2.3}. Let $\ce(-1,0)$ be the dual
bundle to $\ce(1,0)$ and define the \textit{Fefferman space} $\tilde
M$ of $M$ to be space of real lines in $\ce(-1,0)$. We state this more
precisely as follows. Let $\Cal F$ be the bundle obtained by removing
the zero section in $\ce(-1,0)$. Since $\ce(-1,0)$ is a complex line
bundle, we get a free right action of $\Bbb C^*$ on $\Cal F$ which is
transitive on each fibre. Restricting this action to the subgroup
$\Bbb R^*$, we can define $\tilde M$ as the quotient $\Cal F/\Bbb
R^*$. Hence $\tilde M\to M$ is a principal fibre bundle with structure
group $\Bbb C^*/\Bbb R^*\cong U(1)$.

\begin{thm*}
  Let $\tilde M$ be the Fefferman space of $(M,H,\ce(1,0))$. Then the
  CR Cartan bundle $\Cal G\to M$ can be naturally viewed as a
  principal bundle over $\tilde M$ with structure group $G\cap\tilde
  P$.  The normal CR Cartan connection $\om\in\Om^1(\Cal G,\frak g)$
  also gives a Cartan connection on $\Cal G\to\tilde M$.

Denoting by $(p,q)$ the signature of $(M,H)$, the Fefferman space
$\tilde M$ naturally carries a conformal spin structure of signature
$(2p+1,2q+1)$.
\end{thm*}
\begin{proof}
  We have noted in \ref{2.3} that we may identify $\ce(-1,0)$ with
  $\Cal G\x_P\Bbb V^1$. By construction, we can therefore
  view $\tilde M$ as the associated bundle $\Cal G\x_P\Cal P_{\Bbb
    R}\Bbb V^1$ with fibre the space of real lines in $\Bbb V^1$.
  Since $G$ acts transitively on the cone of nonzero null vectors, $P$
  acts transitively on the space of real lines in $\Bbb V^1$. By
  definition (and as observed in \ref{2.2} above) the stabiliser of
  one of these lines is $G\cap\tilde P$, whence $\Cal P_{\Bbb R}\Bbb
  V^1\cong P/(G\cap\tilde P)$. Now $\Cal G\x_P(P/(G\cap\tilde P))$ can
  be naturally identified with the orbit space $\Cal G/(G\cap\tilde
  P)$. Hence we can view $\Cal G$ as a principal bundle over $\tilde
  M$ with structure group $G\cap\tilde P$.
  
  Let $\om\in\Om^1(\Cal G,\fg)$ be the normal CR Cartan connection. By
  definition, $\om$ provides a trivialisation of $T\Cal G$, which is
  equivariant for the action of the structure group $P$ and reproduces
  the generators (in $\frak p$) of fundamental vector fields. But
  $P$--equivariancy of course implies equivariancy for the actions of
  the subgroup $G\cap\tilde P\subset P$, and the fundamental vector
  fields on $\Cal G\to\tilde M$ are exactly those fundamental vector
  fields on $\Cal G\to M$ whose generators lie in $\frak g\cap\tp$.
  Hence $\om$ also defines a Cartan connection on $\Cal G\to\tilde M$.
   From \ref{2.2} we know that the inclusion $G\hookrightarrow\tilde G$
  induces a diffeomorphism $G/(G\cap\tilde P)\to \tilde G/\tilde P$.
  Denoting by $\tilde{\frak g}$ the Lie algebra of $\tilde G$, we
  obtain an inclusion $\frak g\to\tilde{\frak g}$ which induces a
  linear isomorphism $\frak g/(\frak g\cap\tilde{\frak p})\to
  \tilde{\frak g}/\tilde{\frak p}$.  By construction, this isomorphism
  is equivariant for the actions of $G\cap\tilde P$ on both sides
  coming from the adjoint action. The Cartan connection $\om$ on $\Cal
  G\to M$ gives rise to an identification of the tangent bundle
  $T\tilde M$ with the associated bundle $\Cal G\x_{G\cap\tilde
    P}\frak g/(\frak g\cap\tilde{\frak p})$. Since $\tilde G\cong
  SO(2p+2,2q+2)$ and $\tilde P$ is the stabiliser of a null line in
  the standard representation, it is well--known that the natural
  action of $\tilde G$ on $\tilde G/\tilde P$ respects an oriented
  conformal structure of signature $(2p+1,2q+1)$. This gives rise to
  an orientation and a conformal class of inner products of the same
  signature on $T_{e\tilde P}(\tilde G/\tilde P)=\tilde{\frak
    g}/\tilde{\frak p}$, which are invariant under the natural action
  of $\tilde P$. Via the above isomorphism, we obtain corresponding
  structures on $\frak g/(\frak g\cap\tilde{\frak p})$ which are
  invariant under $G\cap\tilde P$.  Passing to the associated bundle,
  we obtain an oriented conformal structure on $\tilde M$.
  
  Thus it remains to construct a natural spin structure. We only
  discuss this for $p,q>0$, if $p=0$ or $q=0$, the argument is
  similar. The group $\tilde G\cong SO_0(2p+2,2q+2)$ is homotopy
  equivalent to its maximal compact subgroup $(O(2p+2)\x O(2q+2))\cap
  SO_0(2p+2,2q+2)$.  Hence the fundamental group of $\tilde G$ is
  $\Bbb Z_2\x\Bbb Z_2$ and the spin group $Spin(2p+2,2q+2)$ is the
  two--fold covering corresponding to the diagonal subgroup. On the
  other hand, the group $G\cong SU(p+1,q+1)$ is homotopy equivalent to
  its maximal compact subgroup, which is isomorphic to $S(U(p+1)\x
  U(q+1))$. Each of the unitary groups has fundamental group $\Bbb Z$
  and hence $\pi_1(G)=\{(k,-k):k\in\Bbb Z\}\subset\Bbb Z\x\Bbb Z$.
  Then the homomorphism $\pi_1(G)\to\pi_1(\tilde G)$ induced by the
  inclusion $G\hookrightarrow\tilde G$ is reduction modulo two in both
  components. Hence its image lies in the diagonal subgroup which
  means that there is a lift to an inclusion $G\hookrightarrow
  Spin(2p+2,2q+2)$.
  
  Restricting this lift to $G\cap\tilde P$ we obtain a homomorphism to
  the stabiliser $\tilde P^{Sp}\subset Spin(2p+2,2q+2)$ of the chosen
  real null line in $\Bbb V$. It is well known that the two--fold
  covering map from this stabiliser onto $\tilde P$ projects onto the
  covering $Spin(\tg/\tp)\to SO_0(\tg/\tp)$. Hence denoting by $K$ the
  kernel of the composition $G\cap\tilde P\to\tilde P^{Sp}\to
  Spin(\tg/\tp)$, we conclude that $\Cal G/K\to \tilde M$ defines a
  two fold covering of a subbundle of the conformal frame bundle of
  $\tilde M$. By construction this is compatible with the projection
  $Spin(\tg/\tp)\to SO_0(\tg/\tp)$, and hence defines a conformal spin
  structure on $\tilde M$.
\end{proof}

\subsection*{Remark}
(1) Although the theorem holds without assuming integrability, we will
soon restrict to the case of CR structures, for which there is a nice
relationship between the normal Cartan connections associated to the
two structures. In the integrable case, the theorem is a minor
variation of known results. The main difference between our
construction and the ones in \cite{BDS}, \cite{LeeF},
\cite{Kuranishi2} is the use of the additional bundle $\ce(1,0)$
instead of the anticanonical bundle.  Thus we obtain an $n+2$--fold
covering of the spaces constructed in those articles. An immediate
advantage of this is the existence of a canonical spin structure as
proved above.

\noindent
(2) Another advantage of our construction is that it exactly recovers
Fefferman's original construction from \cite{Fefferman} for boundaries
of strictly pseudoconvex domains. This can be proved directly by
showing, similarly to the conformal case treated in \cite{confamb},
that the CR standard tractor bundle together with its filtration,
Hermitian metric and connection, can be constructed from the ambient
tangent bundle, the ambient (pseudo--K\"ahler) metric, and its
Levi--Civita connection. Since this needs several non--trivial
verifications, it will be taken up elsewhere.

\noindent
(3) We will give an explicit formula for the conformal structure on
$\tilde M$ in Proposition \ref{4.8} below.

\subsection{Standard tractors}\label{2.5} 
Having extended the structure group of the CR Cartan bundle to $P$, we
can work with the standard representation $\Bbb V$ of $G$ (which does
not make sense for $G/Z(G)$). This leads to standard tractors which
will provide a simple relationship between the CR structure on $M$ and
the conformal structure on $\tilde M$, and at the same time keep track
of the associated Cartan connections. See \cite{TAMS} for the general
relation between tractor bundles and Cartan geometries and
\cite{confamb} for more details on the conformal case.

Restricting the standard representation of $G$ to the subgroup $P$, we
obtain the associated bundle $\Cal T:=\Cal G\x_P\Bbb V$. This is
called the \textit{standard tractor bundle} of $(M,H,\ce(1,0))$. By
construction, it is a rank $n+2$ complex vector bundle over $M$, which
comes equipped with a Hermitian inner product $h$ of signature
$(p+1,q+1)$ that is induced by $\langle\ ,\ \rangle$. The
$P$--invariant subspace $\Bbb V^1\subset\Bbb V$ gives rise to a
natural subbundle $\Cal T^1\subset\Cal T$, which is a complex line
bundle isomorphic to $\Cal E(-1,0)$. Moreover, the fibres of $\Cal
T^1$ are all null with respect to $h$. Since $P\subset SU(\Bbb V)$, a
choice of a nonzero element $\tau\in\La^{n+2}\Bbb V$ induces a
trivialisation of the highest complex exterior power $\La^{n+2}\Cal
T$.

Of course, we can restrict the standard representation further to the
subgroup $G\cap\tilde P\subset P$ and obtain an associated vector
bundle $\tct:=\Cal G\x_{G\cap\tilde P}\Bbb V\to\tilde M$. The
Hermitian inner product on $\Bbb V$ is $G$--invariant, so it gives
rise to a Hermitian bundle metric on $\tct$ of signature $(p+1,q+1)$.
Taking the real part of this defines a real bundle metric $\tilde h$
of signature $(2p+2,2q+2)$ on $\tct$. The real line $\Bbb V_{\Bbb
  R}^1\subset\Bbb V$, stabilised by $G\cap\tilde P$, gives rise to a
real line subbundle $\tct^1\subset\tct$ and each of these lines is
null with respect to $\tilde h$. Thus, defining $\tct^0$ to be the
real orthogonal complement of $\tct^1$, we obtain a filtration
$\tct=\tct^{-1}\supset\tct^0\supset\tct^1$ by smooth subbundles. The
real volume form $\tau\wedge\overline{\tau}$ on $\Bbb V$ induces a
trivialisation of the highest real exterior power $\La^{2n+4}\tct$.

\begin{thm*}
  The Cartan connection $\om$ on $\Cal G$ induces a tractor connection
  $\nabla^{\tct}$ on the bundle $\tct\to \tM$, and $(\tct,\tct^1,\tilde
  h,\nabla^\tct)$ is a standard tractor bundle for the natural
  conformal structure on $\tilde M$. The tractor connection $\nabla^\tct$ is
  normal if and only if the almost CR structure $(M,H)$ is integrable.
\end{thm*}
\begin{proof}
  As a representation of $\tilde P$, the $(2n+2)$nd tensor power of
  $\Bbb V^1_{\Bbb R}$ is dual to the highest exterior power
  $\La^{2n+2}\tg/\tp$. This immediately implies that $\tct^1$ is a
  density bundle on $\tilde M$, which (in the conventions of
  \cite{confamb} but using a tilde to indicate density bundles on
  $\tilde M$) is $\tce[-1]$. In the same way we conclude that
  $\tct/\tct^0\cong\tce[1]$ and $\tct^0/\tct^1\cong T\tilde
  M\otimes\tce[-1]$. Finally the relation between the $\tilde
  P$--invariant conformal class of inner products on $\tg/\tp$ and the
  standard representation shows that the conformal class on $\tilde M$
  constructed in theorem \ref{2.4} comes from the metric on $T\tilde
  M\otimes\tce[-1]\cong\tct^0/\tct^1$ induced by $\tilde h$.
  
  Since the representation $\Bbb V$ of $G\cap\tilde P$ is the
  restriction of a representation of $G$, we may invoke the mechanism
  of \cite[Theorem 2.7]{TAMS} to get a tractor connection on the
  associated bundle $\tct$ from the Cartan connection
  $\om\in\Om^1(\Cal G,\frak g)$. For a smooth section $s\in\Ga(\tct)$,
  consider the corresponding $G\cap\tilde P$--equivariant smooth
  function $f:\Cal G\to\Bbb V$. For a smooth vector field $\xi$ on
  $\tilde M$, choose a lift $\bar\xi\in\frak X(\Cal G)$, and define
  the function corresponding to $\nabla^{\tct}_\xi s$ by $u\mapsto
  (\bar\xi\cdot f)(u)+\om(\bar\xi(u))(f(u))$, where in the second
  summand we use the action of $\frak g$ on $\Bbb V$. As in the proof
  of \cite[theorem 2.7]{TAMS} one concludes that this defines a linear
  connection on $\tct$. By construction, this connection is compatible
  with the bundle metric $\tilde h$ and the volume form on $\tct$.
  Sections of $\tct^1$ correspond to functions with values in $\Bbb
  V^1_{\Bbb R}$. If $f$ is such a function, then so is $\bar\xi\cdot
  f$ for any vector field $\bar\xi$. If $s$ is the corresponding
  section and $s(x)\neq 0$, we see that $\nabla^{\tct}_\xi
  s(x)\in\tct^1_x$ if and only if $\om(\bar\xi)\in\frak g\cap\tp$, the
  stabiliser in $\frak g$ of $\Bbb V^1_{\Bbb R}$. But this means that
  $\bar\xi$ projects to zero on $\tilde M$, so $\xi(x)=0$. Hence, we
  have verified that $\nabla^{\tct}$ is non--degenerate in the sense
  of \cite[2.5]{TAMS}, compare also with \cite[2.2]{confamb}, and thus
  it is a tractor connection. Hence we see that $(\tct,\tct^1,\tilde
  h,\nabla^{\tct})$ is a standard tractor bundle for the natural
  conformal class on $\tilde M$ in the sense of \cite[2.2]{confamb}.
  
  To discuss normality, we need to compare the curvatures of
  $\nabla^{\Cal T}$ and $\nabla^{\tct}$. This is best done in the
  picture of the curvature function.  Initially, the curvature of
  $\om$ is defined as the $\frak g$--valued two form on $\Cal G$ given
  by $(\xi,\eta)\mapsto d\om(\xi,\eta)+[\om(\xi),\om(\eta)]$. This is
  easily seen to be horizontal and $P$--equivariant, so it can be
  interpreted as a two form $\ka$ on $M$ with values in $\Cal
  G\x_P\frak g=\frak{su}(\Cal T)$.  The Cartan connection $\om$ gives
  rise to an isomorphism $TM\cong\Cal G\x_P(\frak g/\frak p)$.
  Therefore, $\ka$ may be viewed as a section of the associated bundle $\Cal
  G\x_P(\La^2(\frak g/\frak p)^*\otimes\frak g)$. The curvature
  function is the $P$--equivariant function $\Cal G\to \La^2(\frak
  g/\frak p)^*\otimes\frak g$ corresponding to this section. We also
  write $\ka$ for the curvature function.
  
  Likewise, $T\tilde M\cong\Cal G\x_{G\cap P}(\tg/\tp)$, so the
  corresponding curvature function has values in $\La^2(\frak g/(\frak
  g\cap\tp))^*\otimes\tg$. To talk about conformal normality, we just
  have to use the isomorphism $\tg/\tp\cong \frak g/(\frak g\cap\tp)$
  obtained in the proof of Theorem \ref{2.2} to interpret the
  curvature function $\tilde\ka$ as having values in
  $\La^2(\tg/\tp)^*\otimes \tg$.
  
  By \cite[Proposition 2.9]{TAMS}, the curvature of the tractor
  connection induced by a Cartan connection is induced by the
  curvature of the Cartan connection. Hence both $\ka$ and $\tilde\ka$
  are induced by the curvature of $\om$, which implies that for
  $u\in\Cal G$, the map $\tilde\ka(u):\La^2(\tg/\tp)\to \tg$ is the
  composition
  \begin{equation}
    \label{curvatures}
    \La^2(\tg/\tp)\overset{\cong}{\longrightarrow} \La^2\big(\frak g/(\frak
    g\cap\tp)\big)\to \La^2(\frak g/\frak p )\overset{\ka(u)}{\longrightarrow} 
\frak g\hookrightarrow\tg.
  \end{equation}
  Now it is well known that a normal conformal Cartan connection must
  be torsion free, which means that if $\nabla^\tct$ is normal, then
  $\tilde\ka(u)$ must have values in $\tp\subset\tg$. This is only
  possible if $\ka(u)$ has values in $\frak g\cap\tp\subset\frak
  p\subset\frak g$. This means that the Cartan connection $\om$ has to
  be torsion free, which implies that $(M,H)$ is integrable, see
  \cite[4.16]{CS}.

To prove the other implication, first note that there is an abelian
subalgebra $\tp_+\subset\tp$, which consists of all maps that
annihilate $\Bbb V^1_{\Bbb R}$ and map the real orthocomplement $\Bbb
V^0_{\Bbb R}$ of this line to $\Bbb V^1_{\Bbb R}$. It is easy to see
that $\tp_+$ is the annihilator of $\tp$ with respect to the trace
form. Thus by the non-degeneracy of the trace form it gives rise to an
isomorphism $\tp_+\cong (\tg/\tp)^*$.  Likewise, there is a subalgebra
$\frak p_+\subset\frak p$, which consists of all maps that annihilate
the complex line $\Bbb V^1$ and map its complex orthocomplement $\Bbb
V^0$ to $\Bbb V^1$. The (real) trace form on $\frak g$ induces an
isomorphism $\frak p_+\cong(\frak g/\frak p)^*$.

Now the conformal normalisation condition can be stated as follows.
Take dual bases $\{\tilde X_j\}$ of $\tg/\tp$ and $\{\tilde Z_j\}$ of
$\tp_+$. Then for each $u\in\Cal G$ and each $X\in\tg/\tp$ the
expression
\begin{equation}\label{confnorm}
\textstyle\sum_{j=0}^{2n+1}[\tilde Z_j,\tilde\ka(u)(\tilde X,\tilde X_j)]
\end{equation}
has to vanish. To obtain appropriate bases, we choose elements
$X_0,\dots,X_{2n}\in\frak g\subset \tg$ which project onto a basis of $\frak
g/\frak p$, 
and put $\tilde X_j:=X_j+\tp$.
Next, take $\tilde X_{2n+1}:=i\id+\tp$.  Under the isomorphism
$\tg/\tp\to\frak g/(\frak g\cap\tp)$ the element $\tilde X_{2n+1}$
corresponds to a nonzero element of $\frak p/(\frak g\cap\tp)$, so
$\{\tilde X_j:j=0,\dots,2n+1\}$ is a basis of $\tg/\tp$. Let $\{\tilde
Z_j\}$ be the dual basis of $\tp_+$. As a linear map on $\Bbb V$ we
can decompose each of the $\tilde Z_j$ uniquely as $Z_j+\hat Z_j$,
where $Z_j$ is complex linear and $\hat Z_j$ is conjugate linear.

First observe that by \eqref{curvatures}, $\tilde\ka(u)$ factors
through $\La^2(\tg/\tp)\to\La^2(\frak g/\frak p)$, which, since
$X_{2n+1}\in\frak p$, implies that the term with $j=2n+1$ does not
contribute to \eqref{confnorm}. For the same reason, it suffices to
take in (\ref{confnorm}) $\tilde X=X+\tp$ for $X\in\frak g$. Further,
we know that $\ka(u)$ has values in $\frak g$, and in particular is
complex linear. Therefore,
$$
[\tilde Z_j,\tilde\ka(u)(\tilde X,\tilde X_j)]=
[Z_j,\tilde\ka(u)(\tilde X,\tilde X_j)]+
[\hat Z_j,\tilde\ka(u)(\tilde X,\tilde X_j)]
$$
is the decomposition into complex linear and conjugate linear parts.
Using \eqref{curvatures} we conclude that the complex linear part of
\eqref{confnorm} is given by
\begin{equation}\label{crnorm1}
\textstyle\sum_{j=0}^{2n}[Z_j,\ka(u)(X+\frak p,X_j+\frak p)].
\end{equation}
Now since each $X_j$ is complex linear, the real trace of $X_j\o\hat
Z_k$ vanishes for all $k$, hence we obtain the real traces tr$(X_j\o
Z_k)=\delta_{jk}$. For $j=2n+1$, we have $X_j=i\id$, which shows that
for $k\leq 2n$ the map $Z_k$ has vanishing complex trace, so it lies
in $\frak{su}(\Bbb V)$. Take a nonzero element $v\in\Bbb V^1_{\Bbb
  R}$.  Then $iv$ lies in the real orthocomplement of $\Bbb V^1_{\Bbb
  R}$, so by the definition of $\tp_+$, $\tilde Z_k(iv)=av$ for some
$a\in\Bbb R$. But then $X_{2n+1}\o\tilde Z_k$ maps $iv$ to $aiv$, and
looking at an appropriate basis, that extends $\{v,iv\}$, one sees
that $\operatorname{tr}(X_{2n+1}\o\tilde Z_k)=2a$, so $a=0$ for $k\leq
2n$. Hence $\tilde Z_k$ vanishes on the complex line $\Bbb V^1$, so
the same is true for $Z_k$. Since $(\Bbb V^1)^\perp$ is a complex
subspace of $\Bbb V$ which is contained in the real orthocomplement of
$\Bbb V^1_{\Bbb R}$, it is mapped to $\Bbb V^1_{\Bbb R}$ by $\tilde
Z_k$ and hence to $\Bbb V^1$ by $Z_k$ for $k\leq 2n$. Hence we have
verified that $\{Z_0,\dots,Z_{2n}\}\subset\frak p_+$, and this is a
basis, which is dual to the basis $\{X_j+\frak p:j\leq 2n\}$ of $\frak
g/\frak p$.

But if $\ka$ is the curvature function of a torsion free normal
parabolic geometry of type $(G,P)$, then \nn{crnorm1} always vanishes.
This is shown in the proof of Theorem 3.8 of \cite{chains}.
Alternatively, it follows by translating the last part of the proof of
Theorem \ref{4.2} below (which is independent of the current
considerations) into a statement on the curvature function. The same
arguments show that $\ka(u)$ has values in $\frak g\cap\tp\subset\frak
p$.  This shows that each summand in \nn{confnorm} lies in
$[\tp_+,\tp]=\tp_+$, and we have just seen that the complex linear
part of the whole sum vanishes. Hence the linear map $\ph:\Bbb
V\to\Bbb V$ defined by \nn{confnorm} is conjugate linear and contained
in $\tp_+$, and we claim that this already implies $\ph=0$.  

The complex subspace $(\Bbb V^1)^\perp$ of $\Bbb V$ is contained in
the real orthocomplement of $\Bbb V^1_{\Bbb R}$. Since $\ph\in\tp_+$,
this implies that $\ph((\Bbb V^1)^\perp)\subset \Bbb V^1_{\Bbb R}$.
But by conjugate linearity, $\ph((\Bbb V^1)^\perp)$ is a complex
subspace of $\Bbb V$, so $\ph$ vanishes on $(\Bbb V^1)^\perp$. Hence
fixing a nonzero element $v\in \Bbb V^1_{\Bbb R}$, there must be an
element $w\in\Bbb V$ such that $\ph(x)=\langle v,x\rangle w$. If $x$
is chosen in such a way that $\langle v,x\rangle=i$, then
$\ph(x)\in\Bbb V^1_{\Bbb R}$, so $w=iav$ for some $a\in\Bbb R$. But
now one immediately verifies that $x\mapsto\langle v,x\rangle iav$ is
symmetric for $\langle\ ,\ \rangle_{\Bbb R}$, so $\ph\in\frak{so}(\Bbb
V)$ only if $\ph=0$.
\end{proof}

\subsection*{Remark}
(1) The standard tractor bundle $(\tct,\tct^1,\tilde h,\nabla^\tct)$
is equivalent to a principal bundle $\tcg\to\tilde M$ with structure
group $\tilde P$ endowed with a Cartan connection
$\tilde\om\in\Om^1(\tcg,\tg)$. The bundle $\tilde{\Cal G}$ can be
constructed as an adapted frame bundle for $\tct$, see
\cite[2.2]{confamb}, which implies that $\tcg=\Cal G\x_{G\cap\tilde
  P}\tilde P$. Since this implies $\tct=\tcg\x_{\tilde P}\Bbb V$, we
obtain a Cartan connection $\tilde\om\in\Om^1(\tcg,\tg)$ by
\cite[Theorem 2.7]{TAMS}. By construction, $\tilde\om$ restricts to
$\om$ on $T\Cal G\subset T\tcg|_{\Cal G}$, which uniquely determines
$\tilde\om$ by the defining properties of Cartan connections.

\noindent
(2) Without the assumption of integrability, the tractor connection
    $\nabla^\tct$ differs from the conformal normal standard tractor
    connection. While we will restrict to the integrable case for the
    rest of this paper, there is scope to use our
    results in the non--integrable case. The route to this should be
    to compute the difference to the normal tractor connection
    explicitly, say in terms of the Nijenhuis tensor of $(M,H)$, and
    then translate the results below (most of which do hold for the
    induced tractor connection in the non--integrable case) to results
    for the normal tractor connection.

\section{Conformal geometry of Fefferman spaces}\label{3}

\subsection{The canonical complex structure on standard
  tractors}\label{3.1} Let $(M,H,\ce(1,0))$ be a CR manifold with
Fefferman space $\tilde M$.  Then we have defined the conformal
standard tractor bundle $\tct$ on $\tilde M$ as $\Cal G\x_{G\cap\tilde
  P}\Bbb V$. Consider multiplication by $i$ as a linear map from
$\Bbb V$ to itself. Since $\langle\ ,\ \rangle$ is Hermitian, this map
lies in $\frak s\frak o(\Bbb V,\langle\ ,\ \rangle_{\Bbb R})=\tg$.
Moreover, the action of any element of $G\cap\tilde P$ commutes with
this map since $G$ consists of complex linear maps. Consequently, the
corresponding constant map $\Cal G\to\tg$ is $(G\cap\tilde
P)$--equivariant and hence gives rise to a section $\Bbb J$ of the
associated bundle $\Cal G\x_{G\cap\tilde P}\tg\cong\tcg\x_{\tilde
  P}\tg$, the \textit{adjoint tractor bundle} $\tca$ of $\tilde M$.
Observe that by construction $\tca=\frak{so}(\tct)$, so the standard
tractor connection induces a linear connection $\nabla^{\tca}$ on
$\tca$, called the \textit{adjoint tractor connection}. Since the
tangent bundle $T\tilde M$ is the associated bundle $\tcg\x_{\tilde
  P}\tg/\tp$, there is a natural projection $\tilde\Pi:\tca\to T\tilde
M$.

\begin{thm*}
  Let $M$ be a CR manifold with Fefferman space $\tilde M$, and let
  $\Bbb J\in\Ga(\tca)$ be the section constructed above. Then we have:

\noindent
(1) $\Bbb J$ makes $\tct$ into a complex vector bundle, it is
orthogonal for $\tilde h$, and $\nabla^{\tca}\Bbb J=0$.

\noindent
(2) The vector field $\k:=\tilde\Pi(\Bbb J)\in\frak X(\tilde M)$ is
nowhere vanishing and generates the vertical bundle of $\tilde M\to
M$. For the conformal Cartan curvature $\tilde\ka\in\Om^2(M,\tca)$ we
have $i_\k\tilde\ka=0$ and $\k$ is a conformal Killing field.
\end{thm*}
\begin{proof}
  (1) Since $\Bbb J$ corresponds to multiplication by $i$ in $\Bbb V$,
  it clearly satisfies $\Bbb J^2=-\id$. Since $\tilde h$ corresponds
  to a Hermitian form on $\Bbb V$, $\Bbb J$ is orthogonal (or
  equivalently skew symmetric). By the definition of $\Bbb J$, if
  $s\in\Ga(\tct)$ corresponds to $f:\Cal G\to\Bbb V$, then $\Bbb Js$
  corresponds to $if$. Now for any tangent vector $\xi$ on $\Cal G$, we
  have $\xi\cdot(if)=i(\xi\cdot f)$ and $\om(\xi)$ is complex linear.
  By definition of $\nabla^{\tct}$ this implies that
  $\nabla^{\tct}\Bbb Js=\Bbb J\nabla^{\tct} s$ for any section $s$.
  Since $\nabla^{\tca}$ is induced by $\nabla^{\tct}$ this shows that
  $\nabla^{\tca}\Bbb J=0$.

\noindent
(2) Fix an element $F\in\frak p$ which acts by multiplication by $i$
on $\Bbb V^1$. Then $i\id-F$ acts trivially on $\Bbb V^1$ and thus
lies in $\tp$. Consequently, the isomorphism $\frak g/(\frak
g\cap\tp)\to\tg/\tp$ induced by the inclusion $\frak g\hookrightarrow
\tg$ maps $F+(\fg\cap\tp)$ to $i\id+\tp$. Since $F\in\frak p$ but
$F\notin (\fg\cap\tp)$, this implies that $i\id+\tp$ is a nonzero
element in the kernel of the projection $\tg/\tp\to \frak g/\frak p$,
which represents $T\pi:T\tilde M\to TM$. Hence the vector field $\k$
is nowhere vanishing and therefore generates the vertical bundle of
$\tilde M\to M$. We have already observed in the proof of Theorem
\ref{2.5} that $\tilde\ka$ comes from the tractor curvature on $\Cal
T$. Since $\k$ lies in the vertical subbundle of $\tilde M\to M$, this
implies $i_{\k}\tilde\ka=0$. Hence $\Bbb J$ satisfies
\begin{equation} \label{cK}
\nabla^{\tca}\Bbb J+i_{\tilde\Pi(\Bbb J)}\ka=0 
\end{equation}
and this is equivalent to $\tilde\Pi(\Bbb J)$
    being a conformal Killing field, compare with
   \cite[Proposition 2.2]{Gopowers}, \cite[Proposition 3.2]{deformations}. 
\end{proof}

This result has some immediate consequences. Using the tractor metric,
we can identify the bundle $\tca=\frak{so}(\tct)$ with the real second
exterior power $\La^2\tct$. Since $\Bbb J$ is a complex structure
then, obviously, as a section of $\La^2\tct$ it is non--degenerate,
i.e.~the $(n+2)$--fold wedge product of $\Bbb J$ with itself is a
nowhere vanishing section of the real line bundle $\La^{2n+4}\tct$.
Hence for each $1\leq k<n+2$ we obtain a nonzero section $\Bbb
J\wedge\dots\wedge\Bbb J$ ($k$ factors) of the bundle $\La^{2k}\tct$.
Since the normal tractor connections on the exterior powers of $\tct$
are induced by $\nabla^{\tct}$, all these sections are parallel.

The filtration $\tct^1\subset\tct^0\subset\tct$ from \ref{2.5} induces
a filtration of the exterior powers of $\tct$ (see e.g.
\cite{BrGodeRham}). Generalising the projection
$\tct\to\tct/\tct^0\cong\ce[1]$, there is a natural projection
$\La^{j}\tct\to\La^{j-1}T^*\tilde M\otimes\ce[j]$. Hence the parallel
section $\Bbb J\wedge\dots\wedge\Bbb J$ gives rise to a weighted
$(2k-1)$--form on $\tilde M$. There is a conformally invariant first
order differential operator defined on sections of $\La^{j-1}T^*\tilde
M\otimes\ce[j]$, which is called the conformal Killing operator, since
for $j=1$ its solutions are conformal Killing fields, see
\cite{Semmelmann} and references therein.

The conformal Killing operator factors through the
composition of the induced connection on $\La^j\tct$ with an invariant 
differential operator which splits the projection
$\La^{j}\tct\to\La^{j-1}T^*\tilde M\otimes\ce[j]$. Moreover, any
parallel section of $\La^j\tct$ is obtained by applying the
splitting operator to its projection. In particular, parallel sections
correspond to special solutions of the conformal Killing equation,
which are called normal conformal Killing forms in \cite{Leitner}. All
these facts are an extremely special case of the machinery of BGG
sequences, whose general version has been developed in \cite{BGG} and
\cite{David-Tammo}. 

\begin{cor*}
Let $\tilde M$ be a Fefferman space. For $1\leq k\leq 2n+1$ let $A_k$
be the space of normal conformal Killing $k$--forms on $\tilde
M$. Then $A_k\neq\{0\}$ for odd $k$ and there is a natural map $A_k\to
A_{k+2}$, which is injective for $k<n+1$ and surjective for $k>n+1$. 
\end{cor*}
\begin{proof}
  From above we know that $A_k$ is isomorphic to the space of parallel
  sections of $\La^{k+1}\tct$. Hence $\Bbb J\wedge\dots\wedge\Bbb J$
  ($k$ factors) projects to a nonzero element of $A_{2k-1}$. Taking
  the wedge product with $\Bbb J$, maps parallel sections of
  $\La^j\tct$ to parallel sections of $\La^{j+2}\tct$, and hence
  induces a map $A_j\to A_{j+2}$. The injectivity and surjectivity
  properties are purely algebraic consequences of the non--degeneracy
  of $\Bbb J$.
\end{proof}
An explicit formula for the odd degree conformal Killing forms which
exist on any Fefferman space is given in Corollary \ref{4.4} below.

\subsection{Relating tractor bundles}\label{3.2}
We next discuss the relation between sections of natural vector
bundles on $M$ and on $\tilde M$. Natural vector bundles on
$(M,H,\ce(1,0))$ are in bijective correspondence with representations
of $P$ via forming associated bundles to the Cartan bundle. Similarly,
natural vector bundles on a conformal manifold are determined by
representations of $\tilde P$. Given a representation of $\tilde P$ on
a vector space $W$, then, in the special case of a Fefferman space, we
have $\tcg\x_{\tilde P} W\cong\Cal G\x_{G\cap\tilde P}W$, so it is
only the restriction of the representation to $G\cap\tilde P$ that
matters.

Assume that we have given representations of $P$ and of $\tilde P$ on
a vector space $W$, which are compatible in the sense that their
restrictions to $G\cap\tilde P$ coincide. Then sections of $\Cal G\x_P
W\to M$ are in bijective correspondence with $P$--equivariant
functions $\Cal G\to W$. On the other hand, sections of
$\tcg\x_{\tilde P}W\cong\Cal G\x_{G\cap\tilde P}W$ are in bijective
correspondence with $(G\cap\tilde P)$--equivariant functions $\Cal
G\to W$. Hence $\Ga(\Cal G\x_P W\to M)$ may be identified with a
subspace of $\Ga(\tcg\x_{\tilde P}W\to\tilde M)$.

There is a simple source for compatible representations: Since $\tilde
G$ contains both $P$ and $\tilde P$ as subgroups, we may use the
restrictions of representations of $\tilde G$ to the two subgroups.
Given such a representation on a vector space $\Bbb W$, the associated
vector bundle $\tcw=\tcg\x_{\tilde P}\Bbb W\to\tilde M$ is the
(conformal) $\Bbb W$--tractor bundle.  Similarly, since $\Bbb W$ is
also a representation of $G$ by restriction, the downstairs bundle
$\Cal W=\Cal G\x_P\Bbb W\to M$ is a (CR) tractor bundle on $M$.

In this way, we obtain a simple relation between conformal tractor
bundles on $\tilde M$ and CR tractor bundles on $\tilde M$. This is
the central tool developed in this article. Of course historically,
and for many applications, irreducible bundles are the geometric
objects studied. However, in contrast to the situation with tractors,
the relation between irreducible bundles on $M$ and on $\tilde{M}$ is
typically complicated. In many cases of interest, such relations can
be deduced from the tractor picture.

The relation between the two types of tractor bundles can be made
explicit using the canonical normal tractor connections. These are
linear connections induced by the canonical Cartan connections
$\tilde\om$ and $\om$, which exist on each tractor bundle, see
\cite{TAMS}. We have noted in Remark \ref{2.5} that $\tilde\om$ is
determined by the fact that its restriction to $T\Cal G$ coincides
with $\om$. As in the case of the standard representation discussed in
\ref{2.5}, this implies that, viewing $\tcw$ as $\Cal G\x_{G\cap\tilde
  P}\Bbb W$, the connection $\nabla^{\tcw}$ is obtained from $\om$ via
the mechanism of \cite[2.7]{TAMS}.

\begin{prop*}
Let $(M,H,\ce(1,0))$ be a CR manifold with Fefferman space $\tilde M$,
$\Bbb W$ a representation of $\tilde G$, and $\tcw\to\tilde M$ and
$\Cal W\to M$ the corresponding tractor bundles. Let $\k=\tilde\Pi(\Bbb
J)\in\frak X(\tilde M)$ be the vector field constructed in \ref{3.1}. Then we
have:

\noindent
(1) A section $\ph\in\Ga(\tcw)$ lies in $\Ga(\Cal W)$ if and only if
    $\nabla^{\tcw}_{\k}\ph=0$. 

\noindent
(2) The restriction of $\nabla^{\tcw}$ to $\Ga(\Cal
    W)\subset\Ga(\tcw)$ descends to a linear connection on the bundle
    $\Cal W\to M$, which coincides with the normal tractor connection
    $\nabla^{\Cal W}$. 
\end{prop*}
\begin{proof}
(1) Since $\tcw\cong\Cal G\x_{G\cap\tilde P}\Bbb W$, sections of
    $\tcw$ are in bijective correspondence with $(G\cap\tilde
    P)$--equivariant functions $\Cal G\to\Bbb W$. A section $\ph$
    lies in the subspace $\Ga(\Cal W)$ if and only if the
    corresponding function $f$ is actually $P$--equivariant. Since
    $P/(G\cap\tilde P)$ is connected, this equivariancy can be checked
    infinitesimally. It is equivalent to the fact that for each $u\in
    \Cal G$ and one (or equivalently any) element $A\in\frak p\setminus
    (\fg\cap\tp)$ we have $\ze_A(u)\cdot f=-A(f(u))$. In the
    left hand side of this expression the fundamental vector field
    differentiates $f$, while on the right hand side $A$ acts by the
    infinitesimal representation on $f(u)\in\Bbb W$. Since $A\in\frak
    p$ we have $\ze_A(u)=\om_u^{-1}(A)$. Choosing for $A$ the element
    $F$ from the proof of Theorem \ref{3.1} we see from the definition
    of the tractor connection that $\om_u^{-1}(F)\cdot f+F (f(u))$
    is exactly  the value at $u$ of the function representing
    $\nabla_{\k}\ph$. 

\noindent
(2) Let $\xi\in\frak X(M)$ be a vector field and let $\tilde\xi$ be a
lift of $\xi$ to $\tilde M$. For $\ph\in\Ga(\Cal W)\subset\Ga(\tcw)$
we have $\nabla^{\tcw}_{\k}\ph=0$, so since $\k$ spans the vertical
bundle of $\tilde M\to M$, we see that $\nabla^{\tcw}_{\tilde\xi}\ph$
depends only on $\xi$ and not on the choice of the lift. From the fact
that $\tilde\xi$ is projectable, one immediately concludes that
$[\k,\tilde\xi]$ lies in the vertical subbundle of $\tilde M\to M$, so
$\nabla^{\tcw}_{[\k,\tilde\xi]}\ph=0$. The curvature $R^{\tcw}$ of
$\nabla^{\tcw}$ is induced by the Cartan curvature of $\tilde\om$.
Hence $R^{\tcw}(\k,\tilde\xi)(\ph)=0$ by part (2) of Theorem
\ref{3.1}, and expanding the definition of the curvature we conclude
that $\nabla^{\tcw}_{\k}\nabla^{\tcw}_{\tilde\xi}\ph=0$.  Therefore,
$\nabla^{\tcw}_{\tilde\xi}\ph$ is an element of $\Ga(\Cal W)$ which we
denote by $\nabla_\xi\ph$. It is straightforward to verify that this
defines a linear connection on $\Cal W\to M$ which by construction
coincides with the tractor connection induced by $\om$.
\end{proof}

In particular, we can apply this result to the standard representation
$\Bbb V$ to characterise sections of $\Cal T\to M$ among sections of
$\tct\to\tilde M$. Note that, by construction for $s,t\in\Ga(\Cal
T)\subset\Ga(\tct)$ the function $\tilde h(s,t)$ is constant along the
fibres of $\pi:\tilde M\to M$ and descends to the real part of
$h(s,t)$. 

\subsection{Twistor spinors on Fefferman spaces}\label{3.3}
We next describe a surprising application of Proposition \ref{3.2}. It
is known that a certain variant of Fefferman spaces always admits
nontrivial twistor spinors. This has been shown for CR dimension one
in \cite{Lewandowski} and in general in \cite{Baum} by direct and
involved computations. Here we obtain a conceptual proof, without
any computations, via a simple analog of the proof of the existence of
parallel spinors on Calabi--Yau manifolds. 

Suppose that $\Bbb W$ is a representation of $\tilde G$ as in
\ref{3.2} above. Even if $\Bbb W$ is irreducible as a representation
of $\tilde G$ it may well happen that as a representation of $G$,
$\Bbb W$ decomposes into several irreducible components. If $\Bbb
W=\Bbb W_1\oplus\dots\oplus\Bbb W_\ell$ as a representation of $G$,
then this decomposition is also valid as a representation of
$G\cap\tilde P$. Hence for a Fefferman space, the corresponding
conformal tractor bundle $\tcw$ admits a decomposition into a direct
sum of bundles. Each of the summands corresponds to a CR tractor
bundle $\Cal W_i\to M$. For each $i$, we can view $\Ga(\Cal W_i)$ as a
subspace of $\Ga(\tcw)$.  Since the tractor connection on $\tcw$ can
be obtained from the Cartan connection $\om$ on $\Cal G$, this
decomposition is compatible with the tractor connection, and the
restriction of $\nabla^\tcw$ to the subspace $\Ga(\Cal W_i)$ descends
to the normal CR tractor connection on $\Cal W_i\to M$.

\begin{cor*}
Let $(M,H,\ce(1,0))$ be a CR manifold with Fefferman space $\tilde
M$. Then, for the canonical spin structure on $\tilde M$, there is a two
parameter family of nonzero twistor spinors. 
\end{cor*}
\begin{proof}
  From the proof of Theorem \ref{2.4} we may take $\tilde G$ to be the
  spin group. Denoting by $\Bbb S$ the basic spin representation, the
  resulting conformal tractor bundle $\tcs\to\tilde M$ is known as
  the bundle of local twistors, which was introduced by Penrose (see
  \cite{Penrose} and references therein) in dimension four. The
  corresponding tractor connection is the local twistor transport and
  it is well--known that parallel sections of this bundle are in
  bijective correspondence with twistor spinors \cite{Fried89,BFGK}.  
Hence it suffices to
  show that the bundle $\tcs\to\tilde M$ admits a two--parameter
  family of parallel sections. But the restriction of the basic spin
  representation to $SU(p+1,q+1)\subset Spin(2p+2,2q+2)$ splits into a
  direct sum of irreducibles among which there are two copies of the
  trivial representation (see
  e.g.~\cite{Wang}). Hence we obtain two copies of
  $C^\infty(M,\Bbb C)$ in $\Ga(\tcs)$, on which the spin tractor
  connection descends to the exterior derivative.  Hence constant
  functions in these two copies give rise to a two parameter family of
  parallel sections of $\tcs$.
\end{proof}

\subsection{Relating adjoint tractors}\label{3.4}
An important example of a conformal tractor bundle, which is in general
indecomposable but which admits  a direct sum splitting on Fefferman spaces,
is the adjoint tractor bundle $\tca$. To apply the ideas of \ref{3.3}
we have to understand the restriction of the adjoint representation of
$\tilde G=SO_0(2p+2,2q+2)$ to the subgroup $G=SU(p+1,q+1)$. Given a
real linear map $f:\Bbb V\to\Bbb V$ which is skew symmetric with
respect to $\langle\ ,\ \rangle_{\Bbb R}$, we can uniquely split $f$
into a complex linear part $f_1$ and a conjugate linear part. Then
$f_1$ is skew Hermitian with respect to $\langle\ ,\ \rangle$ so it
can be written as the sum of an element of $\frak s\frak u(\Bbb V)$
and a real multiple of $i\id$. On the other hand, mapping $\ph$ to
$\langle\ ,\ph(\ )\rangle$ defines a linear isomorphism between the
space of those conjugate linear endomorphisms of $\Bbb V$ which are also skew
symmetric with respect to $\langle\ ,\ \rangle_{\Bbb R}$, and the space
of complex bilinear skew symmetric maps $\Bbb V\x\Bbb V\to\Bbb C$.
Hence we see that $\tg=\frak g\oplus\Bbb R\oplus\La^2_{\Bbb C}\Bbb
V^*$ and this decomposition is invariant under the action of $G$. Since
all the summands are irreducible for $G$, this is the complete
decomposition.

From this elementary representation theory, it follows that the
adjoint tractor bundle $\tca$ splits into three pieces that are
preserved by $\nabla^{\tca}$. It is easy to make this splitting
explicit. Let $\{\ ,\ \}$ be the algebraic bracket on $\tca$ induced
by the commutator of endomorphisms of $\tct$, and let $B$ be the real
trace form, i.e.~the non--degenerate bilinear form mapping two
endomorphisms to the real trace of their composition. By construction
we then have $B(\Bbb J,\Bbb J)=-(2n+4)$. Given a section
$s\in\Ga(\tca)$ we can write the conjugate linear part of $s$ as
$\tfrac{1}{4}\{\Bbb J,\{s,\Bbb J\}\}$. Finally the trace part of $s$
only sits in the complex linear part and is given by
$\tfrac{-1}{2n+4}B(s,\Bbb J)\Bbb J$. Hence we conclude that the full
decomposition of $s\in\Ga(\tca)$ is given as
\begin{equation}
s=\left(s-\tfrac{1}{4}\{\Bbb J,\{s,\Bbb J\}\}+
\tfrac{1}{2n+4}B(s,\Bbb J)\Bbb J\right)+
\tfrac{-1}{2n+4}B(s,\Bbb J)\Bbb J+\tfrac{1}{4}\{\Bbb J,\{s,\Bbb J\}\}~. 
\label{tca-split}\end{equation}
Since $\nabla^\tca$ satisfies a Leibniz rule with respect to $\{\ 
,\ \}$ and $\Bbb J$ and $B$ are parallel, we see that, as expected,
this decomposition is preserved by $\nabla^{\tca}$. In particular, we
see that for the CR adjoint tractor bundle $\Cal A\to M$ we have
$$
\Ga(\Cal A)=\{s\in\Ga(\tca):\{s,\Bbb J\}=0,B(s,\Bbb
J)=0,\nabla^{\tct}_{\k}s=0\} 
$$
and the restriction of $\nabla^\tca$ to this subspace descends to the
normal adjoint tractor connection on $\Cal A$. 

We can apply this to obtain a complete description of infinitesimal
conformal isometries of a Fefferman space.  As noted in the proof of
Theorem \ref{3.1}, infinitesimal conformal isometries of $\tilde M$
are in bijective correspondence with smooth sections $s\in\Ga(\tca)$
such that $\nabla^{\tca}s+i_{\tilde\Pi(s)}\tilde\ka=0$.
First we need the following result on $\tilde{M}$.

\begin{lem*}
  Suppose that $s\in\Ga(\tca)$ satisfies
  $\nabla^{\tca}s+i_{\tilde\Pi(s)}\tilde\ka=0$. Then $\{s,\Bbb J\}$ is
  a parallel section of $\tca$, and the function $B(s,\Bbb J)$ is
  constant. 
\end{lem*}
\begin{proof}
  For $\xi\in\frak X(M)$ we obtain using the Leibniz rule and that
  $\Bbb J$ is parallel
$$
\nabla^\tca_\xi\{s,\Bbb J\}=\{\nabla^\tca_\xi s,\Bbb
J\}=-\{\tilde\ka(\Pi(s),\xi),\Bbb J\},
$$
which vanishes since $\tilde\ka$ has complex linear values. 

Similarly, naturality of $B$ implies that 
$$
\xi\cdot B(s,\Bbb J)=B(\nabla^\tca_\xi s,\Bbb
J)=-B(\tilde\ka(\Pi(s),\xi),\Bbb J),
$$
which vanishes since $\tilde\ka$ has values in $\frak{su}(\tct)$. 
\end{proof}

\begin{thm*}
  Let $(M,H,\ce(1,0))$ be a CR manifold with Fefferman space $\tilde
  M$. Then the space of conformal Killing fields on $\tilde M$
  naturally splits into a direct sum $A_1\oplus\Bbb R\k\oplus A_2$,
  where $A_1$ is isomorphic to the space of infinitesimal CR
  automorphisms of $M$. The space $A_1$ can be identified with the
  space of solutions of a second order CR invariant linear
  differential operator defined on the bundle $Q=TM/H$. The space
  $A_2$ can be identified with a subspace of the joint kernel of two
  CR invariant linear first order differential operators defined on
  $H\otimes\ce(-1,1)$.
\end{thm*}
\begin{proof}
  For a conformal Killing field consider the corresponding section
  $s\in\Ga(\tca)$ and its decomposition $s=s_1+a\Bbb J+s_2$ according
  to \eqref{tca-split}. By the Lemma, the function $a$ is constant and
  $\nabla^\tca s_2=0$. By  \cite[Proposition 2.2]{Gopowers} or 
\cite[Corollary 3.5]{deformations}, this
  implies $i_{\tilde\Pi(s_2)}\tilde\ka=0$ and therefore
  $i_{\tilde\Pi(s)}\tilde\ka=i_{\tilde\Pi(s_1)}\tilde\ka$. Hence $s_1$
  and $s_2$ independently satisfy the infinitesimal automorphism
  equation, and $\tilde\Pi(s)=\tilde\Pi(s_1)+a\k+\tilde\Pi(s_2)$ is
  the decomposition claimed in the theorem. 
  
  Since $\k$ hooks trivially into $\tilde\ka$ we see that
  $\nabla^{\tca}_\k s_1=0$, so $s_1\in\Ga(\Cal A)$. Since the
  connection $\nabla^{\tca}$ and the curvature $\tilde\ka$ descend to
  their CR counterparts, $s_1$ satisfies $\nabla^{\Cal
    A}s_1+i_{\Pi(s_1)}\ka=0$, where $\Pi:\Cal A\to TM$ is the natural
  projection and $\ka$ is the CR tractor curvature. This equation is
  equivalent to $s_1$ giving rise to an infinitesimal CR automorphism,
  see \cite[Proposition 3.2]{deformations}. The interpretation of
  infinitesimal CR automorphisms as solutions of an invariant operator
  follows using the BGG machinery, see \cite[Theorem
  3.4]{deformations}.
  
 Since $s_2$ is conjugate linear it may be viewed as lying in
  $\Ga(\La^2_{\Bbb C}\Cal T^*)$ and it is parallel for the normal
  tractor connection on that bundle. In analogy with the exterior
  powers of the conformal standard tractor bundle in \ref{3.1}, there
  is a natural projection on $\La^2_{\Bbb C}\Cal T^*$ to an
  irreducible quotient bundle, which here is $H\otimes\ce(-1,1)$. The
  machinery of BGG sequences can be used to construct invariant
  differential operators between irreducible bundles from tractor
  connections. The first operators in a BGG sequence are defined on
  the irreducible quotient of the tractor bundle in question. The
  other bundles occurring in the sequence are irreducible subquotients
  of the bundles of forms with values in the initial tractor bundle.
  The isomorphism type of these subquotients can be computed using
  representation theory. If we start from the tractor bundle
  $\La^2_{\Bbb C}\Cal T^*$ then these computations show that one
  obtains two invariant operators defined on $H\otimes\ce(-1,1)$, one
  having values in $S^2H\otimes\ce(-2,0)$ and the other having values
  in the tensor product of tracefree endomorphisms of $H$ with
  $\ce(-1,1)$. Again by representation theory, these two invariant
  operators must be of first order. 
  
  A particular consequence of the construction of BGG sequences is
  that projecting a parallel section of a tractor bundle to the
  irreducible quotient, one always obtains a section which lies in the
  kernel of the first operators in the BGG sequence. Hence $s_2$
  projects onto a section of $H\otimes\ce(-1,1)$ which is annihilated
  by the two invariant operators discussed above.
\end{proof}

\subsection*{Remark}
(1) The decomposition of conformal Killing fields is described explicitly
in \ref{4.8b} below. 

\noindent
(2) There is a general classification of first order invariant
    differential operators on arbitrary parabolic geometries, see
    \cite{Slovak-Soucek}. In particular, this implies that the two
    first order operators occurring in the theorem are both given by
    taking one Webster--Tanaka derivative and then projecting to the
    given irreducible component. 

\subsection{Relating densities and weighted tractors}\label{3.5}
Another source of representations of $P$ and $\tilde P$ which are
compatible in the sense of \ref{3.2} is provided by density bundles.
Let $\rho:P\to\Bbb C^*$ be the representation defined by the action of
$P$ on $\Bbb V^1$. Then we can form the representation $g\mapsto
\rho(g)^{-w}\overline{\rho(g)}^{-w'}$, and the corresponding
associated bundle is the density bundle $\ce(w,w')\to M$. Restricting
this representation to $G\cap\tilde P$ we obtain multiplication by
$\la^{-w-w'}$ (for $\la\in {\Bbb R}^*$), so the corresponding
associated bundle is the complexified density bundle $\tce_{\Bbb
  C}[w+w']:=\tce[w+w']\otimes\Bbb C$.

Recall that for any natural vector bundle on a conformal manifold, one
has the \textit{fundamental derivative} or fundamental $D$--operator,
see \cite[section 3]{TAMS}. We will denote this operator by
$\tilde{\Bbb D}$. In particular we get an operator, on $\tilde{M}$,
$\Ga(\tca)\x\Ga(\tce_{\Bbb C}[w])\to\Ga(\tce_{\Bbb C}[w])$ written as
$(s,\al)\mapsto\tilde{\Bbb D}_s\al$.

\begin{prop*}
A complex density $\al\in\Ga(\tce_{\Bbb C}[w+w']\to\tilde M)$ lies in
the subspace $\Ga(\ce(w,w')\to M)$ if and only if $\tilde{\Bbb D}_{\Bbb
J}\al=(w-w')i\al$.
\end{prop*}
\begin{proof}
Consider the $(G\cap\tilde P)$--equivariant function $f:\Cal G\to\Bbb
C$ representing $\al$. As in the proof of Proposition \ref{3.2} we
conclude that $\al\in\ce(w,w')$ if and only if $\om_u^{-1}(F)\cdot
f=-(w'-w)if(u)$, where $F\in\frak p$ is the element from the proof of
Theorem \ref{3.1}. Of course, the right--hand side represents the
function $(w-w')i\al$, so it suffices to verify that the left--hand
side represents $\tilde{\Bbb D}_{\Bbb J}\al$.

To see this, let us view $\Cal G$ as a subset of $\tcg$ and extend $f$
equivariantly to a function $\tilde f:\tcg\to\Bbb C$. By definition of
$\tilde{\Bbb D}$, the density $\tilde{\Bbb D}_{\Bbb J}\al$ is
represented by the function $(\tilde\om^{-1}\o \ph)\cdot\tilde f$,
where $\ph:\tcg\to\tg$ is the function representing $\Bbb J$. For
$u\in\Cal G$ we have $\ph(u)=i\id$ and hence this function evaluates
in $u$ to $\tilde\om_u^{-1}(i\id)\cdot\tilde f$. As we have noted in
the proof of Theorem \ref{3.1}, we have $F-i\id\in\tp$. Since this
element annihilates the real line $\Bbb V^1_{\Bbb R}$, it acts
trivially on powers of this, viz.\ real one dimensional
representations.  Since $\tilde f$ is $\tilde P$--equivariant, we
conclude that $\tilde\om^{-1}(F-i\id)\cdot\tilde f=0$ and hence
$\tilde\om^{-1}(i\id)\cdot\tilde f=\tilde\om^{-1}(F)\cdot\tilde f$.
Since $\tilde\om$ restricts to $\om$ on $T\Cal G\subset T\tcg|_{\Cal
  G}$ we conclude that $\tilde\om_u^{-1}(F)=\om_u^{-1}(F)$ for
$u\in\Cal G$, and since this is tangent to $\Cal G$ we obtain
$\om_u^{-1}(F)\cdot\tilde f=\om_u^{-1}(F)\cdot f$.
\end{proof}

Using this result we can now easily extend the characterisation of
Proposition \ref{3.2} to weighted tractor bundles, i.e.~tensor
products of tractor bundles with density bundles. Given a
representation $\Bbb W$ of $\tilde G$, we consider the tractor bundles
$\Cal W\to M$ and $\tcw\to\tilde M$ as in \ref{3.2}. We define $\Cal
W(w,w'):=\Cal W\otimes\ce(w,w')$ and $\tcw_{\Bbb
  C}[w+w']:=\tcw\otimes\ce_{\Bbb C}[w+w']$. If $\Bbb W$ is a complex
representation then the tensor product can be taken over $\Bbb C$,
otherwise it is understood as a real tensor product. The basic
operator on such bundles is the \textit{double--D--operator}
\cite{GoSrni99,Goadv,TAMS}, which is obtained by coupling the
fundamental derivative to the tractor connection. This is well
defined, since both the fundamental derivative and the tractor
connection satisfy a Leibniz rule. In particular, the conformal
double--D defines an operator $\tilde{\Bbb
  D}^\nabla:\Ga(\tca)\otimes\Ga(\tcw_{\Bbb C}[w+w'])\to\Ga(\tcw_{\Bbb
  C}[w+w'])$. Explicitly, for $s\in\Ga(\tca)$, $\al\in\Ga(\ce_{\Bbb
  C}[w+w'])$ and $\ph\in\Ga(\tcw)$ we have
$$
\tilde{\Bbb D}^\nabla_s(\al\otimes\ph)=(\tilde{\Bbb
  D}_s\al)\otimes\ph+\al\otimes\nabla^{\tcw}_{\Pi(s)}\ph. 
$$
The above with Proposition \ref{3.2} gives the following result.
\begin{cor*}
A section $\ph\in\Ga(\tcw_{\Bbb C}[w+w'])$ lies in the subspace
$\Ga(\Cal W(w,w'))$ if and only if $\tilde{\Bbb D}^\nabla_{\Bbb
  J}\ph=(w-w')i\ph$. 
\end{cor*}

\subsection{Relating fundamental derivatives and
  double--D--operators}\label{3.6} The fundamental derivative is
defined on any natural bundle on a manifold endowed with an arbitrary
parabolic geometry. The CR version of these operators can be coupled
to CR tractor connections to obtain a CR version of
double--D--operators. Both types of operators play a central role in
the invariant theory and invariant operator theory of conformal, CR
geometry \cite{GoSrni99,Goadv} (and other parabolic geometries
\cite{TAMS}). Our next task is to relate the CR versions of these
operators to the conformal variants on the Fefferman space.  The basis
for this is a simple identity in the conformal setting.
\begin{lem*}
Let $\tcw\to\tilde M$ be any conformally natural bundle on a Fefferman
space, and let $\Bbb J\in\Ga(\tca)$ be the canonical complex
structure. Then for arbitrary sections $s\in\Ga(\tca)$ and
$\ph\in\Ga(\tcw)$ we get
$$
\tilde{\Bbb D}_\Bbb J\tilde{\Bbb D}_s\ph=\tilde{\Bbb D}_s\tilde{\Bbb
  D}_\Bbb J\ph+\tilde{\Bbb D}_{(\nabla^{\tca}_\k s-\{\Bbb
  J,s\})}\ph. 
$$
\end{lem*}
\begin{proof}
  By Section 3.7 of \cite{TAMS}, the difference $\tilde{\Bbb D}_\Bbb
  J\tilde{\Bbb D}_s\ph-\tilde{\Bbb D}_s\tilde{\Bbb D}_\Bbb J\ph$ is
  given by $\tilde{\Bbb D}_{[\Bbb J,s]}\ph$, where $[\Bbb J,s]$
    denotes the Lie bracket of the two adjoint tractor fields. A
    formula for this Lie bracket is given in \cite[Proposition
    3.6]{TAMS}, and together with Proposition 3.2 of that reference
  this gives
$$
[\Bbb J,s]=\tilde{\Bbb D}_{\Bbb J}s-\nabla^{\tca}_{\Pi(s)}\Bbb
J-\tilde\ka(\Pi(\Bbb J),\Pi(s))=\tilde{\Bbb D}_{\Bbb J}s,
$$
where we have used that $\Bbb J$ is parallel and $\Pi(\Bbb J)=\k$
hooks trivially into the Cartan curvature. Applying again Proposition
\ref{3.2} of \cite{TAMS} we get $\tilde{\Bbb D}_{\Bbb
  J}s=\nabla^{\tca}_\k s-\{\Bbb J,s\}$, which implies the result.
\end{proof}

From this, we immediately get the relation between
double--D--operators.

\begin{prop*}
  Consider a representation $\Bbb W$ of $\tilde G$ and the
  corresponding weighted tractor bundles on $M$ and $\tilde M$ as in
  \ref{3.5}. Then for any $s\in\Ga(\Cal A)\subset\Ga(\tca)$ the
  double--D--operator $\tilde{\Bbb D}^\nabla_s$ acting on
  $\Ga(\tcw_{\Bbb C}[w+w'])$ preserves the subspace $\Ga(\Cal
  W(w,w'))$ and by restriction gives the CR double--D--operator $\Bbb
  D^\nabla_s$ on that subspace.
\end{prop*}
\begin{proof}
  The description of $\Ga(\Cal A)\subset\Ga(\tca)$ in \ref{3.4} shows
  that $\tilde\nabla^\tca_\k s=\{\Bbb J,s\}=0$. By the lemma, the
  fundamental derivatives with respect to $\Bbb J$ and to $s$ commute.
  By the proof of the lemma, this also implies that $[\Bbb J,s]=0$ and
  using Proposition 3.6 of \cite{TAMS} this gives $[\k,\Pi(s)]=0$.
  Since $\k$ hooks trivially into $\tilde\ka$, and the curvature of
  any tractor connection is given by the action of $\tilde\ka$, we
  conclude that also $\nabla^{\tcw}_\k$ and $\nabla^{\tcw}_{\Pi(s)}$
  commute.  Hence $\tilde{\Bbb D}^\nabla_s$ commutes with $\tilde{\Bbb
    D}^\nabla_{\Bbb J}$, and this with Corollary \ref{3.5} shows that
  the subspace $\Ga(\tcw_{\Bbb C}[w+w'])$ is preserved by $\tilde{\Bbb
    D}^\nabla_s$.
  
  In \ref{3.2} we have already seen that conformal tractor connections
  restrict to their CR counterparts. Hence to see that $\tilde{\Bbb
    D}^\nabla_s$ restricts to its CR counterpart it suffices to show
  that $\tilde{\Bbb D}_s$ restricts to $\Bbb D_s$ on
  $\Ga(\ce(w,w'))\subset\Ga(\tce_{\Bbb C}[w+w'])$. The fact that
  $s\in\Ga(\Cal A)\subset\Ga(\tca)$ implies that the corresponding
  equivariant function $f:\tcg\to\tg$ has the property that $f(\Cal
  G)\subset {\frak g}\subset \tg$. Consequently, the vector field
  $\xi\in \frak X(\tcg)$ characterised by $\tilde\om(\xi(u))=f(u)$ for
  all $u\in\tcg$ has the property that its restriction to $\Cal G$ is
  tangent to $\Cal G$. Since $\tilde\om$ restricts to $\om$ on vectors
  tangent to $\Cal G$, we see that $\xi|_{\Cal G}\in\frak X(\Cal G)$
  is the vector field associated to $s\in\Ga(\Cal A)$ via the CR
  Cartan connection. Now the result follows immediately from the
  definition of the fundamental derivative.
\end{proof}

\subsection{Complexified adjoint tractors}\label{3.7}
A deeper understanding and richer theory of fundamental derivatives
and double--$D$--operators is revealed, in this context, by passing to
the complexification $\tca_{\Bbb C}$ of the adjoint tractor bundle.
Recall from \ref{3.4} that $\tg=\frak g\oplus\Bbb R\oplus\La^2_{\Bbb
  C}\Bbb V^*$ as a $\frak g$--module, with the first two summands
corresponding to complex linear maps and the last summand
corresponding to conjugate linear maps. Since the last summand is
already complex, its complexification splits into a holomorphic and an
anti--holomorphic part. The holomorphic part is isomorphic to
$\La^2_{\Bbb C}\Bbb V^*$ and since $\Bbb V^*\cong\overline{\Bbb V}$
via the Hermitian form, the anti--holomorphic part is isomorphic to
$\La^2_{\Bbb C}\Bbb V$. Moreover, the isomorphism between conjugate
linear maps $\ph$ in $\tg$ and skew symmetric complex bilinear maps on
$\Bbb V$ is given by $\ph\mapsto \langle \ ,\ph(\ )\rangle$, which
implies that the complex structure on $\La^2_{\Bbb C}\Bbb V^*$
corresponds to mapping $\ph$ to $v\mapsto -i\ph(v)$.

Passing to associated bundles, we see that the subspace of elements of
$\tca_{\Bbb C}$ which anti--commute with $\Bbb J$ is isomorphic to
$\La^2_{\Bbb C}\tct^*\oplus\La^2_{\Bbb C}\tct$. The two summands are
characterised by $is=-\Bbb J\o s$ respectively $is=\Bbb J\o s$. Note
that $is=\mp\Bbb J\o s$ implies that $\{\Bbb J,s\}=\mp 2is$. 

Since we deal with unweighted tractor bundles here, the
double--$D$--operator $\tilde{\Bbb D}^\nabla_{\Bbb J}$ is simply given
by the tractor connection $\tilde\nabla_{\k}$. Now by corollary
\ref{3.5}, a section $s\in\Ga(\La^2_{\Bbb C}\tct^*)$ lies in the
subspace $\Ga(\La^2\Cal T^*(-1,1))\subset\Ga(\La^2_{\Bbb C}\tct^*)$ if
and only if
$$
\tilde\nabla_\k s=\tilde{\Bbb D}^\nabla_{\Bbb J}s=-2is=\{\Bbb J,s\}. 
$$
In the same way, $s\in\Ga(\La^2_{\Bbb C}\tct)$ lies in the subspace
$\Ga(\La^2\Cal T(1,-1))$ if and only if $\tilde\nabla_\k s=\{\Bbb
J,s\}$. 

Now consider a complex weighted tractor bundle $\tcw_{\Bbb C}[w+w']$
as in \ref{3.6}. For a section $\ph$ of this bundle, we can view
$\tilde{\Bbb D}^\nabla\ph$ as a section of 
$$
L(\tca,\tcw_{\Bbb C}[w+w'])\cong 
L_{\Bbb C}(\tca_{\Bbb C},\tcw_{\Bbb C}[w+w']).
$$ 
Hence we can form $\tilde{\Bbb D}^\nabla_s\ph$ for all sections
$s\in\Ga(\tca_{\Bbb C})$. In close analogy to the proof of
Proposition \ref{3.6} we then conclude that if $s\in\Ga(\tca_{\Bbb
C})$ sits in either of the subspaces $\Ga(\La^2\Cal T^*(-1,1))$ or
$\Ga(\La^2\Cal T(1,-1))$, then the operator $\tilde{\Bbb D}^\nabla_s$
preserves the subspaces $\Ga(\Cal W(w,w'))\subset \Ga(\tcw_{\Bbb
C}[w+w'])$ for all $w,w'$.

Recall that $\tg\cong\tg^*$, via the Killing form, and the
decomposition $\tg=\frak g\oplus\Bbb R\oplus\La^2\Bbb V^*$ is
orthogonal with respect to the Killing form. These results extends to
the complexification.  In the complexification of the last factor of
the decomposition, the holomorphic and anti--holomorphic parts are
both isotropic with respect to the Killing form, which induces a
duality between the two parts. Using this duality (in the language of
associated bundles), we can interpret $\tilde{\Bbb D}^\nabla$ as an
operator mapping sections of $\tcw_{\Bbb C}[w+w']$ to sections of
$\tca_{\Bbb C}\otimes\tcw_{\Bbb C}[w+w']$, and this target splits
according to the splitting of $\tca_{\Bbb C}$. From Proposition
\ref{3.6} and the above considerations we get:

\begin{thm*}
Let $\tcw_{\Bbb C}[w+w']$ be a weighted complex
tractor bundle. Consider the double--$D$--operator as an operator 
$$
\tilde{\Bbb D}^\nabla:\Ga(\tcw_{\Bbb C}[w+w'])\to \Ga(\tca_{\Bbb
  C}\otimes \tcw_{\Bbb C}[w+w']).
$$
(1) Passing to the tracefree part of the complex linear part in the
$\tca_{\Bbb C}$--component, one obtains an operator which descends to
the CR double--$D$, viewed as 
$$
\Bbb D^\nabla:\Ga(\Cal W(w,w'))\to\Ga(\Cal A\otimes\Cal W(w,w')).
$$ 
(2) If one forms the holomorphic part of the conjugate linear part in
the $\tca_{\Bbb C}$--component, then the result descends to an
operator
$$
\Ga(\Cal W(w,w'))\to\Ga(\La^2\Cal T\otimes\Cal W(w-1,w'+1)). 
$$ 
(3) If one forms the anti--holomorphic part of the conjugate linear
part in the $\tca_{\Bbb C}$--component, then the result descends to an
operator
$$
\Ga(\Cal W(w,w'))\to\Ga(\La^2\Cal T^*\otimes\Cal W(w+1,w'-1)). 
$$
\end{thm*}

The Theorem gives, via the Fefferman space, a geometric interpretation
to the operators in parts (2) and (3) which were constructed directly
(but without a conceptual interpretation) in \cite{GoSrni99}. These
will be described explicitly in \ref{4.9}.

\section{Tractor calculus on a Fefferman space}\label{4}
In this section, we will describe the (complexified) standard tractor
bundle and the tractor calculus on a Fefferman space more explicitly.
To do this, we show that the version of CR tractors introduced in
\cite{Gover-Graham} describes the normal CR standard tractor bundle,
and relate it to the calculus on the Fefferman space.  Hence we obtain
an explicit description in terms of a chosen pseudo--hermitian
structure on the underlying CR manifold.  

\subsection{Pseudo--hermitian structures}\label{4.1}
We review some facts about pseudohermitian structures on a CR manifold
$(M,H)$, primarily to fix the conventions, which follow
\cite{Gover-Graham}.  We will assume that $M$ is orientable, which
implies that the annihilator $H^{\perp}$ of $H$ in $T^*M$ admits a
nonvanishing global section. From the non-degeneracy of the CR
structure such a section $\th$ is a contact form on $M$, and it is
called a \textit{pseudohermitian structure}. We fix an orientation on
$H^{\perp}$ and restrict consideration to $\theta$'s which are
positive relative to this orientation.  The \textit{Levi form} of
$\theta$ is the Hermitian form $\Le$ on $H^{1,0}\subset TM\otimes\Bbb
C$ defined by
$$
\Le(Z,{\ol{W}})=-2id\theta (Z,{\ol{W}}),
$$
so this exactly corresponds to $\Cal L^{\Bbb C}$ introduced in
\ref{2.1} under the trivialisation given by $\th$. 

Given a pseudohermitian structure $\theta$, define the \textit{Reeb
  field} $r$ to be the unique vector field on $M$ satisfying
\begin{equation}\label{Tnorm}
\theta (r)=1 ~~~{\rm and}~~~i_rd\theta=0.
\end{equation}
An \textit{admissible coframe} is a set of complex valued forms $\{
\theta^{\al}\}$, $\al =1,\cdots ,n$, which satisfy $\theta^{\al}(r)=0$
and whose restrictions to $H^{1,0}$ are complex linear and form a
basis for $(H^{1,0})^*$.  We will use lower case Greek indices to
refer to frames for $T^{1,0}$ or its dual.  We may also interpret
these indices abstractly, so will denote by $\ce^{\al}$ the bundle
$H^{1,0}$ (or its space of sections) and by $\ce_{\al}$ its dual, and
similarly for the conjugate bundles or for tensor products thereof.
By integrability and \eqref{Tnorm}, we have
$$
d\theta = ih_{\al\beb} \theta^{\al}\wedge \theta^{\beb}
$$
for a smoothly varying Hermitian matrix $h_{\al\beb}$, which we may
interpret as the matrix of the Levi form in the frame $\theta^{\al}$,
or as the Levi form itself in abstract index notation. Using the
inclusion $Q\hookrightarrow\ce(1,1)$ from \ref{2.3}, the Levi form
$\Cal L^{\Bbb C}$ itself can be viewed as an canonical section of
$\ce_{\al\beb}(1,1)$ which we also denote by $\bh_{\al\beb}$. By
$\bh^{\al\beb}\in\ce^{\al\beb}(-1,-1)$ we denote its inverse.  These
will be used to raise and lower indices without further mention.

By $\nabla$ we denote the \textit{Webster--Tanaka connections} (on
various bundles) associated to $\th$. In particular, these satisfy
$\nabla\th=0$, $\nabla h=0$, $\nd \bh=0$, $\nabla r=0$, and $\nd J=0$,
so the decomposition $T_{\Bbb C}M=H^{1,0}M\oplus H^{0,1}M\oplus\Bbb C
r$ is invariant under $\nabla$.

Therefore, if we decompose a tensor field relative to this splitting
(and/or its dual), we may calculate the covariant derivative
componentwise.  Each of the components may be regarded as a section of
a tensor product of $\ce^{\al}$ or its dual or conjugates thereof.
Therefore we will often restrict consideration to the action of the
connection on $\ce^{\al}$ or $\ce_{\al}$.  We will use indices $\al,
\ol{\al}, 0$ for components with respect to the frame
$\{\theta^{\al},\theta^{\ol{\al}}, \theta \}$ and its dual, so that the
0-components incorporate weights.  If $f$ is a (possibly
density-valued) tensor field, we will denote components of the
(tensorial) iterated covariant derivatives of $f$ in such a frame by
preceding $\nd$'s, e.g.  $\nd_{\al} \nd_0 \cdots \nd_{\beb} f$.  As
usual, such indices may alternately be interpreted abstractly.  So,
for example, if $f_\be\in \ce_\be(w,w')$, we will consider $\nd f$ as
the triple $\nd_{\al}f_\be\in \ce_{\al\be}(w,w')$,
$\nd_{\ol{\al}}f_\be\in \ce_{\ol{\al}\be}(w,w')$, $\nd_0f_\be \in
\ce_\be(w-1,w'-1)$.

\subsection{CR tractor calculus}\label{4.2} 
Our next task is to show that the calculus introduced in
\cite{Gover-Graham} is consistent with the (complexified) CR standard
tractor bundle and connection as described here. The defining feature
of the bundle $\Cal T_{\Ph}$ constructed in \cite{Gover-Graham} is
that any choice of a pseudohermitian structure $\th$ on $M$ gives rise
to an identification 
\begin{equation*}
\Cal T_\Phi \stackrel{\theta}{=} \ce(1,0)\oplus\ce_\al(1,0)\oplus\ce(0,-1).
\end{equation*}
For a section $T_\Phi\in \Cal T_\Phi$ one writes
$$
[T_\Phi]_\theta=\left( \begin{array}{c} \sigma\\
                            \tau_{\be}\\
                 \rho \end{array}\right),
$$
or equivalently 
$$
T_\Phi= \si Y_\Phi + \tau_\be W_\Phi{}^\be + \rho Z_\Phi,
$$
for $\si\in \ce(1,0)$, $\tau_\be \in\ce_\be(1,0)$,
$\rho\in\ce(0,-1)$ and sections $Y_\Phi\in \ct_\Phi(-1,0)$,
$W_\Phi{}^\be \in \ct_\Phi{}^\be(-1,0)$, and $Z_\Phi\in\ct_\Ph(0,1)$
which depend on $\theta$. Changing scale from $\th$ to $\hat\th=e^{\Up}\th$,
the expression for $[T_\Phi]_{\hat\theta}$ is determined by
$$
\wh{W}_\Phi{}^\al = W_\Phi{}^\al +\Up^\al Z_A \quad \quad 
\wh{Y}_\Phi= Y_\Phi-\Up_\be W_\Phi{}^\be  
-\frac{1}{2}(\Up_\be\Up^\be-i\Up_0)Z_\Phi ~,
$$
while $Z_\Phi$ is independent of the choice of $\th$ and e.g.\ $\Up_\al:=\nd_\alpha \Up$. In
particular, this shows that $Z_\Ph$ gives rise to an isomorphism from
$\ce(0,-1)$ onto a subbundle $\Cal T^1\subset\Cal T_{\Ph}$. For two
sections $T_\Ph$ and $T'_\Ph$ the quantity
$\si\overline{\rho'}+\rho\overline{\si'}+\bh^{\al\beb}\tau_\al\tau'_\beb$
is independent of the choice of $\th$ so one obtains a well-defined
hermitian metric $h^{\Ph\overline{\Ps}}$ on $\Cal T_\Ph$. Note that 
 the subbundle $\Cal T^1$ is isotropic for
$h^{\Ph\overline{\Ps}}$. Taking $\Cal T^0$ to be the orthocomplement
of $\Cal T^1$, we obtain $\Cal T^1\subset\Cal T^0\subset\Cal T_{\Ph}$,
and for any choice of $\th$ the elements of $\Cal T^0$ are
characterised by $\si=0$. Projecting onto the $\si$--component shows
that $\Cal T_\Ph/\Cal T^0\cong\ce(1,0)$, while projecting onto the
middle component we have $\Cal T^0/\Cal
T^1\cong\ce_{\be}(1,0)$. The filtration of $\ct_\Ph$ can be
equivalently described as a composition series which we write as 
$$
\ct_\Phi = \ce(1,0)\lpl \ce_\al(1,0)\lpl \ce(0,-1).
$$

Now the bundles $\ce(0,-1)$ and $\ce(1,0)$ are, by definition, conjugate
to $\ce(-1,0)$ respectively $\ce(0,1)$. Via the Levi form, the bundle
$ \ce_{\be}(1,0)$ is identified with the conjugate of
$H^{1,0}\otimes\ce(-1,0)$.  Therefore, the conjugate bundle $\Cal
T_{\overline{\Ph}}$ to $\Cal T_\Ph$, which via $h^{\Ph\overline{\Ps}}$
is identified with the dual bundle $\Cal T^\Ph$, has a composition
series
$$
\ct^\Phi = \ce(0,1) \lpl \ce^\al(-1,0)\lpl \ce(-1,0).
$$
which is exactly as for the standard tractor bundle from \ref{2.5}.
For $\ct_{\overline{\Ph}}$, we obtain a canonical section
$Z_{\overline{\Ph}}\in\ct_{\overline{\Ph}}(1,0)$ which maps to the
line subbundle. Clearly given a choice of contact form we also have
the projectors $Y_{\overline{\Ph}}\in\ct_{\overline{\Ph}}(0,-1)$ and
$W_{\al\overline{\Ph}}\in\ct_{\al\overline{\Ph}}(1,0)$. We can use the
Hermitian metric to raise and lower tractor indices. For example, we
obtain $Z^\Ph\in\ct^{\Ph}(1,0)$ which represents the natural inclusion
$\ce(-1,0)\hookrightarrow\ct^\Ph$ as well as the natural projection
$\ct_{\Ph}\to\ce(1,0)$.

\subsection{Normality}\label{4.3}
The next step in \cite{Gover-Graham} is to introduce a linear
connection on $\Cal T_{\Ph}$. Since we often have to use this tractor
connection coupled to a Webster--Tanaka connection, the best move is
to denote both by $\nabla$. Which connection is acting 
is determined by the objects it acts on, so this should cause no
confusion. In \cite{Gover-Graham} the extension of the linear connection to the
complexified tangent bundle is provided directly: in the display
(3.3) of that reference the authors produce explicit formulae for
$\nabla_\al T_\Ph$, $\nabla_{\beb}T_{\Ph}$ and $\nabla_0T_\Ph$, for a
section $T_\Ph$ in terms of a choice of $\th$, and verify that this
definition is independent of the choice. It is also verified there 
that the connection is Hermitian. On the other hand, by
construction it is compatible with the complex structure on
$\ct_\Ph$.

The connection on $\ct_\Phi$ coupled with the Tanaka-Webster
connection acts on the projectors as follows.
$$
\begin{aligned}
\nd_\be Y_\Phi &= iA_{\al\be} W_\Phi{}^\al+T_\beta Z_\Phi \\ 
\nd_\be W_\Phi{}^\al &= -\delta^\al_\be Y_\Phi - \Rho_\be{}^\al Z_\Phi\\
\nd_\be Z_\Phi &= 0 ~,
\end{aligned}
$$
and
$$
\begin{aligned}
\nd_\beb Y_\Phi &= \Rho_{\al\beb} W_\Phi{}^\al-T_\beb Z_\Phi\\
\nd_\beb W_\Phi{}^\al &= i A_\beb{}^\al Z_\Phi\\
\nd_\beb Z_\Phi &= \bh_{\al \beb} W_\Phi{}^\al ~,
\end{aligned}
$$
and
$$
\begin{aligned}
\nd_0 Y_\Phi &= \frac{i}{n+2}P Y_\Phi+2iT_\al W_\Phi{}^\al+iS Z_\Phi\\
\nd_0 W_\Phi{}^\al &= -i \Rho_\be{}^\al W_\Phi{}^\be +\frac{i}{n+2}PW_\Phi{}^\al 
+2iT^\al Z_\Phi\\
\nd_0 Z_\Phi &= -iY_\Phi+\frac{i}{n+2}\Rho Z_\Phi  ~.
\end{aligned}
$$
where the quantities on the right--hand--side that we have not
defined above are torsion and curvature components of the
Webster--Tanaka connection. The definitions can be found in
\cite{Gover-Graham}.

These formulae then determine the connection of $\ct_\Phi$ in an
obvious way. In particular taking the covariant derivative of $\rho
Z_\Ph$ for a locally nonvanishing section $\rho\in\ce(0,-1)$ and
factoring by $\Cal T^1$, the resulting tensorial map $TM\to \Cal
T_\Ph/\ct^1$ is injective. Indeed passing further to $\Cal
T_\Ph/\ct^0$ exactly extracts the coefficient of the Reeb field by the
formula for $ \nd_0 Z_\Phi$ while the formula for $\nd_\beb Z_\Phi$
shows that the middle component will be injective on $H^{0,1}$.

\begin{thm*}
The bundle $\Cal T_\Ph$ can be naturally identified with the dual of
the normal standard tractor bundle in such a way that the filtration,
the Hermitian metric and the connection $\nabla$ are mapped to their
canonical counterparts. 
\end{thm*}
\begin{proof} Let us write
  $\Cal A\to M$ to denote the bundle of skew Hermitian endomorphisms
  of $\Cal T_\Ph$. The filtration of $\Cal T_\Ph$ gives rise to a
  filtration on $\Cal A$, while the commutator defines a tensorial Lie
  bracket. Hence $\Cal A$ becomes a bundle of filtered Lie algebras
  modelled on $\frak s\frak u(p+1,q+1)$ and thus an abstract adjoint
  tractor bundle for $\Cal T$ in the sense of \cite[Section
  2.2]{TAMS}. The filtration has the form $\Cal A=\Cal
  A^{-2}\supset\Cal A^{-1}\supset\dots\supset\dots\Cal A^2$, and the
  component $\Cal A^j$ is characterised by the facts that for
  $\ell=-1,0,1$ its elements map $\Cal T^i_\Ph$ to $\Cal T^{i+j}_\Ph$,
  where $\Cal T^\ell_\Ph=\Cal T_\Ph$ for $\ell<-1$ and $\Cal
  T^\ell_\Ph=0$ for $\ell>1$. 
  
  A corresponding principal bundle, with structure group the subgroup
  $P\subset SU(p+1,q+1)$ from \ref{2.2}, can be constructed as the
  frame bundle for $\Cal T_\Ph$ sensitive to the filtration structure.
  This is an adapted frame bundle in the sense of \cite[Section
  2.2]{TAMS}.  The connection $\nabla$ on $\Cal T_\Ph$ from above
  induces a connection on $\Cal A$, and one immediately verifies that
  the non--degeneracy property observed above implies that this is a
  tractor connection. In view of Section 2.12 of \cite{TAMS} we
  therefore only have to verify that the curvature $\Om$ of $\nabla$,
  which is computed in \cite{Gover-Graham} satisfies the normalisation
  condition.
  
  An explicit formula for the normalisation condition can be found in
  \cite[2.5]{CS}. Translated to geometric terms, this reads as
\begin{equation}\label{CRnorm}
  0=\sum_{j}\{\eta_j,\Om(\xi,\xi_j)\}+
  \tfrac{1}{2}\sum_{j}\Om(\Pi(\{\eta_j,A\}),\xi_j)
\end{equation}
for all vector fields $\xi$, where $A\in\Cal A$ satisfies
$\Pi(A)=\xi$, the $\xi_j$ form a real local frame for $TM$ and the
$\eta_j$ form the dual frame for $T^*M$. As in \ref{3.4}, the brackets
$\{\ ,\ \}$ denote the tensorial Lie bracket on $\Cal A$ induced by
the Lie bracket on $\frak g$. Moreover, one uses the natural
identification $T^*M\cong\Cal A^1$, see \cite[Section 2.8]{TAMS}.
Since the formula for $\Om$ in \cite{Gover-Graham} refers to a choice
of $\th$, we may assume that $\eta_0=\th$ and $\xi_0=r$ while the
remaining elements form dual frames for $H$ and $H^*$. This implies
that $\eta_0\in\Cal A^2$ and hence $\{\eta_0,A\}\in\Cal A^0=\ker(\Pi)$
for all $A\in\Cal A$. On the other hand, if $\Pi(A)=\xi\in H$, then
$A\in\Cal A^1$, which implies that $\Pi(\{\eta_j,A\})=0$ for all $j$.
For $A\in\Cal A$ such that $\Pi(A)=\xi_0=r$, one computes directly
that $\Pi(\{\eta_j,A\})=-J\xi_j$.  Then vanishing of the second sum in
\eqref{CRnorm} follows from the fact that
$\Om_\al{}^\al{}_\Ph{}^\Ps=0$ which is formula (3.4) of
\cite{Gover-Graham}.

To analyse the first sum in \eqref{CRnorm}, we use the matrix
representation of $\Om$ from \cite{Gover-Graham}. The bracket $\{\ ,\ 
\}$ may then be computed as a commutator of matrices, by representing
a one-form $\ph$ by the matrix 
$$
\begin{pmatrix} 0 & 0
  & 0 \\ \ph_\al & 0 & 0 \\ -i\ph_0 & -\ph_\beb & 0\end{pmatrix}.
$$
From this, and the formula for $\Om$ in \cite{Gover-Graham}, we
immediately conclude that the term with $j=0$ does not contribute to
the sum. Hence we are left with real dual frames of $H$ and $H^*$, and
it suffice to show that the expression vanishes if one sums over
complex dual frames. Then the computation can be done directly, and
the identities for the curvature quantities derived in
\cite{Gover-Graham} immediately show that all traces which show up in
the result vanish.
\end{proof}

\subsection{Conformal standard tractors with a parallel and orthogonal complex
  structure}\label{4.4}

For conventions on conformal structures and results on conformal
tractor calculus, we refer to \cite{BEG,Gover-Peterson}, but we use
tildes in the notation to distinguish conformal objects from CR
objects. Let $\tilde M$ be a smooth manifold of dimension $2n+2$
endowed with a conformal structure $[g]$ of signature
$(2p+1,2q+1)$. By $\bg\in\ce_{(ab)}[-2]$ we denote the conformal
metric. Denoting the standard tractor bundle by $\tct^A$, we have a
conformally invariant metric $\tilde h$ of signature $(2p+2,2q+2)$ and
a composition series
$$
\tct^A=\tce[1]\lpl \tce_a[1]\lpl \tce[-1].
$$
Let $X^A$ be the canonical section of $\tct^A[1]$ which represents the
inclusion $\tce[-1]\to \tct^A$. Tractor indices can be raised and
lowered using $\tilde h$. For example, we obtain a natural section
$X_A\in\tct_A[1]$, which represents the natural projection
$\tct^A\to\tce[1]$. We shall raise and lower indices in this way
without further mention. 

A Weyl structure is a splitting of the filtration of the tractor
bundle. Evidently this is equivalent to a section $Y_A$ of
$\tct_A[-1]$ such that $X^AY_A=1$ and $Y^AY_A=0$. We are most
interested in  splittings that arise from a choice of scale. A
scale is a section $s$ of $\tce_+[1]$, the positive ray subbundle in
$\tce[1]$. (This determines a metric from the conformal class viz.\ 
$g=s^{-2}\bg$.)  There is a conformal generalisation of the exterior
derivative $\tilde{d}$, see \cite{BrGodeRham}. This arises from the
restriction of the exterior derivative on the total space of the
conformal metric bundle to differential forms which are homogeneous
for the obvious ${\Bbb R}_+$-action.  Thus this operator is
conformally invariant, first order and, for example, maps sections of
$\ce[1]$ to sections of $\tct_A/\tce[-1]$. Then $Y_A$ is the unique
null (weighted) tractor which maps to $s^{-1}\tilde d s$ under the
canonical quotient map $\tct_A[-1]\to\tct_A[-1]/\tce[-2]$.  Henceforth
$Y_A$ on $\tilde M$ will mean the section arising from a scale in this
way. Having made this choice, we obtain
$$
\tct^A\stackrel{s}{=} \tce[1]\oplus \tce_a[1]\oplus \tce[-1].
$$
and we write $\bZ$ for the complementary projector/injector.
That is a triple $(\si,\mu_a,\rho)$ from the direct
sum represents the element $\si Y^A +\bZ^A{}_a\mu^a+\rho X^A\in
\tct^A$. Under a change of scale $s\mapsto e^{-\Up} s$ these
projectors transform according to
\begin{equation}\label{XYZtrans}
\textstyle
\begin{array}{rl}
\hat \bZ^{Ab}=\bZ^{Ab}+\Up^bX^A, &
\hat Y^A=Y^A-\Up_b\bZ^{Ab}-\frac12\Up_b\Up^bX^A.
\end{array}
\end{equation}
where $\Up_a=d\Up$. The tractor metric $\tilde h$ is
characterised by $(\si,\mu_a,\rho)\mapsto 2\si \rho+
\bg_{ab}\mu^a\mu^b$. On the other hand, the conformal metric $\bg$ is
recovered from the tractor metric by the expression
\begin{equation}\label{Metric}
\bg_{ab}\xi^a\eta^b=\tilde h_{AB}\bZ^A{}_a\bZ^B{}_b\xi^a\eta^b .
\end{equation}
Note that although the projector $\bZ^A{}_a$ depends on a choice of
metric, from the conformal class, it follows easily from \nn{XYZtrans},
and the inner product relations amongst the projectors, that the
expression on the right-hand side is independent of this choice.  

We use the same symbol $\tilde\nabla$ for the Levi-Civita connections
determined by a choice of scale, and also for the canonical tractor
connections, the distinction is again by context. Using the coupled
connection, the tractor connection is then determined by
\begin{equation}
  \label{conf-connection}
  \tilde{\nabla}_aX^A=\Bbb Z^{A}{}_a \quad \tilde{\nabla}_a\Bbb
  Z^A{}_b=-\widetilde{\Rho}_{ab}X^A-Y^A\bg_{ab}\quad 
  \tilde{\nabla}_aY^A=\widetilde{\Rho}_{ab}\Bbb Z^{Ab}, 
\end{equation}
where $\widetilde{\Rho}_{ab}$ is the conformal Rho--tensor (or
Schouten--tensor).

\newcommand{\J}{{\Bbb J}}
\newcommand{\D}{{\Bbb D}}

Suppose that the standard tractor bundle is endowed with a complex
structure $\Bbb J$ which is orthogonal (or equivalently skew
symmetric) with respect to the tractor metric. In abstract indices we
have $\Bbb J_{A}{}^{B}$ with $\J_A{}^B\J_B{}^C =-\delta_A^C$ and $\Bbb
J_{AB}=-\Bbb J_{BA}$. Using $\Bbb J$ we obtain a canonical section
$\K^A:=(\Bbb JX)^A=X^B\Bbb J_B{}^A\in\tct^A[1]$.  Since $\Bbb J$ is
orthogonal, we immediately obtain $\K^A\K_A=X^AX_A=0$ as well as
$\K^AX_A=-\K^AX_A$ so that $\K$ is null and orthogonal to $X$.

Since $\K^AX_A=0$, the element $\k^a:=\K^A \bZ_A{}^a\in\tce^a$ is
independent of the choice of scale. (In the case of a Fefferman space
this is, by construction, the conformal Killing field $\k$ from
Theorem \ref{3.1}.) Since in any scale $\K^A- \bZ^A{}_b\k^b=X^A \rho$,
for some density $\rho$, it follows from \eqref{Metric}, and that
$\K^A$ is null, that $\k^a$ is null for the conformal structure.

As mentioned $\K$ is, by construction, independent of any choice of
scale.  On the other hand a choice of scale determines a dual object
viz.\ $\L^A:= Y^B\Bbb J_B{}^A$. Arguing in a manner similar to the
above, we see that $\L^A$ is null and orthogonal to $Y$, while
$\L^A\K_A=1$ since $\Bbb J$ is orthogonal.

\subsection*{Definition} A conformal scale $s\in \ce_+[1]$ on $\tilde
M$ is called \textit{preferred} (with respect to the complex structure
$\Bbb J$) if and only if for the fundamental derivative $\tilde{\Bbb
D}$ we have $\tilde{\Bbb D}_{\Bbb J}s=0$.

\begin{prop*} Let $\tilde M$ be a smooth manifold of dimension $2n+2$
  endowed with a conformal structure of signature $(2p+1,2q+1)$ and an
  orthogonal parallel complex structure $\Bbb J$ on the standard
  tractor bundle $\tct$. Let $\K^A$ the canonical weighted tractors
  constructed above and let $\k^a=\K^A\Bbb Z_A{}^a$ be the conformal
  Killing field underlying $\Bbb J$. Let $s$ be a preferred scale on
  $\tilde M$, let $\L^A$ be the associated weighted tractor, put
  $\ell_a:=\L_A\Bbb Z^A{}_a$. Then in the scale $s$ we have:

\noindent
(1) $\tilde\nabla_a\k^a=0$, so $\k^a$ is a Killing field for the
    metric determined by $s$. 

\noindent
(2) $\K^A=\Bbb Z^A{}_b\k^b$, $\L_A=\Bbb Z_A{}^a\ell_a$, and hence
    $\k^a\ell_a=1$ and $\ell^a\ell_a=0$.

\noindent
(3) $\k^b\tilde\nabla_b\k^a=\ell_b\tilde\nabla^b\k^a=0$.

\noindent
(4) $\ell_a=\tilde\Rho_{ab}\k^b$ and $\k^b\tilde\nabla_b\ell_a=0$. 

\noindent
(5) The almost complex structure $\Bbb J$ is explicitly given by 
$$ 
\J_{AB}=2Y_{[A}\Bbb Z_{B]}{}^b\k_b+\Bbb Z_A{}^a\Bbb
  Z_B{}^b\tilde{\nabla}_a\k_b+2X_{[A}\Bbb
  Z_{B]}{}^b\widetilde{\Rho}_{bc}\k^c.
$$
\end{prop*}
\begin{proof}
Recall the definition of the conformally invariant operator
\textit{tractor--$D$} operator, which acts on arbitrary weighted
tractor fields. For a conformal tractor bundle $\tcw$, the operator
$D_A$ maps sections of $\tcw[w]$ to sections of
$\tct_A\otimes\tcw[w-1]$. In our notation, it is given by
\begin{equation}
  \label{tractor-D}
  D_A t:=(2n+2w)wY_At+(2n+2w)\Bbb
  Z_A{}^a\tilde{\nabla}_at-X_A(\tilde{\nabla}^a\tilde{\nabla}_a+
  w\widetilde{\Rho})t,
\end{equation}
where $\widetilde{\Rho}=\widetilde{\Rho}_a{}^a$.

For a vector field $v^a$ define the tractor $V^A\in\tct^A[1]$ as being
given, in a scale $s'$, by $\tilde{V}^A:=\Bbb
Z^A{}_av^a-\frac{1}{2n+2}X^A\tilde{\nabla}'_av^a$, where we use
$\tilde\nabla'$ for covariant derivatives with respect to $s'$.  It is
easily verified that this defines a conformally invariant operation.
From Lemma 2.1 of \cite{Gopowers} the fact that
$\k^a$ is a conformal Killing is equivalent to the corresponding
tractor $K^A$ satisfying $D_A K_B= - D_B K_A$, while from Proposition
2.2 there, the differential splitting operator relating conformal
Killing fields to sections of the adjoint tractor bundle satisfying
equation \nn{cK} is given by $v^a\mapsto
\frac{1}{2n+2}D_{[A}\tilde{V}_{B]}$, where $[\cdots]$ indicates that
we take the skew part over the enclosed indices.  Since $\J$ is
parallel and $\k=\Pi(\J)$, we can recover $\J$ from $\k$ as
$\J^{A}{}_B=\frac{1}{2n+2}D^AK_B$.  But then $\K_B=\J^{A}{}_BX_A=K_B$.
That is in any scale $s'$ we have $\K^A=\Bbb
Z^A{}_a\k^a-\frac{1}{2n+2}X^A\tilde{\nabla}'_a\k^a$.

By the formula for $\tilde{\D}$ in \cite{luminy} the equation
$\tilde{\Bbb D}_{\Bbb J}s=0$ expands to
$$ 
n\k^a\tilde{\nabla}'_a s- s \tilde{\nabla}'_a \k^a=0.
$$ This holds for any $s'$. But, by construction for the metric
determined by $s$, we have $\tilde\nabla_as=0$, and since $s$ is
nowhere vanishing (1) follows.

\newcommand{\h}{{\tilde h}}

Thus, in the preferred scale $s$, we have $\K^A=\bZ^A{}_b\k^b$, and 
$$
0=\K^AY_A= (\Bbb J_B{}^AX^B)Y_A= (\Bbb J_B{}^A Y_A)X^B= -\L_B X^B. 
$$ 
Since we already know that $\L_BY^B=0$, we get $\L_A=\Bbb
Z_A{}^a\ell_a$ and then $\K^A\L_A=1$ and $\L^A\L_A=0$ imply the rest
of (2). The formula for $\J$ in (5) is then obtained by expanding
$\J^{A}{}_B=\frac{1}{2n+2}D^A\K_B$.

Using (5) to expand $\L_A=\Bbb J_B{}^AY^B$, we obtain $\L_A=\Bbb
Z_A{}^b\widetilde{\Rho}_{bc}\k^c$ and hence the formula for $\ell_a$
in (4). Note that this implies $\widetilde{\Rho}_{ab}\k^a\k^b=1$,
which is familiar from Sparling's characterisation of Fefferman
spaces, see \cite{Graham}. 

Differentiating $\k^a\k_a=0$, we get $\k^a\tilde\nabla_b\k_a=0$, which
implies that first equation in (3) by the skew symmetry of
$\tilde\nabla_b\k_a=0$. Next, from the definition of $\L$, we get
$\L^A\Bbb J_{AB}=-Y_B$, and expanding this using (5), the second part
of (3) follows. 

Since $\k^a$ is a Killing field, its Lie derivative annihilates
$\tilde\Rho_{ab}$. Of course the Lie bracket of $\k^a$ with itself vanishes,
and  so the Lie derivative by  $\k^a$ annihilates
 $\ell_a=\tilde\Rho_{ab}\k^b$. This reads as
$0=\k^b\tilde\nabla_b\ell_a- \ell^b\tilde\nabla_b\k_a$, and the second
summand vanishes by (3).
\end{proof}

From the formula in part (5), we immediately get an explicit formula
for the normal conformal Killing forms obtained from $\J$ in Corollary
\ref{3.1}.
\begin{cor*}
  For each $j=1,\dots,n-1$, the form defined in a preferred scale $s$
  as $\k_{[a}(\tilde{\nabla}_{a_1}\k_{b_1})\dots
  (\tilde{\nabla}_{a_j}\k_{b_j]})$ defines a normal conformal Killing
  $(2j+1)$--form on $\tilde M$.
\end{cor*}

\newcommand{\h}{\tilde h}

In a preferred scale the corresponding metric from the conformal class
may be put in the form
\begin{equation}\label{metric}
g_{ab}=2 \k_a\odot \ell_b + \h_{ab}
\end{equation}
where $\h_{ab}$ annihilates $\k^a$ and $\ell^a$. 

In the case of a Fefferman space, we will shortly describe such a
metric in terms of tensors on $M$. Before doing that, we will study
the decomposition of the tangent spaces induced by a choice of
preferred scale in more detail.

\subsection{Decomposition of the tangent bundle}\label{4.4a}
Let us write $V$ for the line subbundle in $T\tM$ spanned by
$\k$. Since $\k$ is null, the orthocomplement $\k^\perp$ contains $V$,
and defining $\tilde H:=\k^\perp/V$, and $\tQ:=T\tilde M/\k^\perp$, we
obtain a composition series for $T\tM$, namely
\begin{equation}\label{Tcomp}
T\tM=\tQ\lpl \tH \lpl V.
\end{equation}
The developments in \ref{4.4} above show that a choice of a preferred
scale $s$ leads to a splitting of the filtration
$V\subset\k^\perp\subset T\tilde M$. Since $\ell_a\k^a=1$, we see that
$\ell^a$ spans a line subbundle in $T\tilde M$ which is complementary
to $\k^\perp$, and that $\tilde H_s:=\k^\perp\cap\ell^\perp$ is a
corank one subbundle of $\k^\perp$ complementary to $V$. In
particular, a choice of preferred scale induces an identification of
$\tilde H_s=\k^\perp\cap\ell^\perp$ with $\tilde H$. Let $\Bbb I^a_b$
be the projector onto $\k^\perp\cap\ell^\perp$, i.e.
\begin{equation}\label{Idef}
\Bbb I^a_b=\delta^a_b-\k^a\ell_b-\ell^a\k_b.
\end{equation}

By duality, the whole picture carries over to the cotangent bundle in
an obvious way. In particular, we obtain a composition series
\begin{equation}\label{ttMcomp}
T^*\tM= V^* \lpl \tH^* \lpl \tQ^* .
\end{equation}
The splitting determined by a preferred scale $s$ comes from the line
subbundle spanned by $\ell_a$ and the annihilator of $\ell^a$. The
corresponding decomposition of a one--form $\om$ explicitly reads as
\begin{equation}\label{formsplit}
\om_a = \om_b\k^b\ell_a + \om_b\ell^b\k_a + \Bbb I^b_a\om_b.
\end{equation}

Combining the projectors $\Bbb I$ and $\Bbb Z$ we obtain 
\begin{equation}\label{tWdef}
\tilde W_A{}^a:=\Bbb Z_A{}^b\Bbb I^a_b=\Bbb Z_A{}^a-\K_A\ell^b-\L_A\k^b.
\end{equation}
Viewed as a map $T^*\tM\to \tct_A$, this annihilates the subbundle
spanned by $\k_a$ and $\ell_b$, and is injective on $\tH^*_s$. In terms
of this the form $\h_{ab}$ in the metric \nn{metric} is given by
$$
\h_{ab}= s^{-2}\tW^A{}_a \tW^B{}_b \tth_{AB}.
$$

\newcommand{\tJ}{{\tilde J}}

Now $\J$ preserves the (weighted) tractor subspace spanned pointwise
by $X$ and $\K$ and similarly (in a choice of preferred scale $s$) it
preserves the subspace spanned pointwise by $Y$ and $\L$. Thus $\J$
determines a canonical complex structure $\tJ$ on the subquotient
bundle $\tH=\k^\perp/V$ of $T\tM$. We may equally view this as an
almost complex structure on $\tH_s=\k^\perp\cap\ell^\perp\subset
T\tilde M$. In this picture, it is given by
$\tJ_a{}^b=\J_A{}^B\tW^A{}_a\tW_B{}^b$. From parts (3) and (4) of
Proposition \ref{4.4} and \eqref{tWdef} we see that $\Bbb Z_A{}^a\Bbb
Z_B{}^b\tilde\nabla_a\k_b=\tW_A{}^a\tW_B{}^b\tilde\nabla_a\k_b$, so
part (5) of Proposition \ref{4.4} immediately gives 
\begin{equation}\label{J2} \J^A{}_B =X^A\L_B - \K^AY_B+\tilde W^A{}_a\tilde
W_B{}^b\tilde{\nabla}^a\k_b + Y^A\K_B- \L^A X_B.
\end{equation}
which in turns shows that $\tJ_a{}^b=\tilde\nabla_a\k^b$.

We will also need the complexified version of the decomposition of the
tangent and cotangent bundles. The composition series \eqref{Tcomp}
and \eqref{ttMcomp} carries over to the complexified setting without
changes. The main additional input is that the complexification of the
subquotient $\tilde H$ splits into a holomorphic and an
anti--holomorphic part. We will use upper Greek indices for the
holomorphic, and bared upper Greek indices for the anti holomorphic
part. Correspondingly, the projectors $\Bbb I^a_b$ give rise to $\Bbb
I^\al_b$ and $\Bbb I^\alb_b$. Correspondingly, we get $\tW_A{}^\al$
and $\tW_A{}^\alb$. 

\newcommand{\bJ}{\Bbb J} 
\newcommand{\tZ}{{\tilde{Z}}}
\newcommand{\Phb}{{\overline{\Phi}}} 
\newcommand{\tY}{{\tilde{Y}}}

\subsection{The complexified standard tractor bundle}\label{4.4b}
To relate the conformal calculus developed so far to CR tractor
calculus, it will convenient to complexify the standard tractor
bundle. This complexification splits into its holomorphic and
anti-holomorphic parts via $t\mapsto (\frac{1}{2}(t-i\J t),
\frac{1}{2}(t+i\J t))$.  We write $\tct^A\otimes\Bbb C=\tct^\Phi\oplus
\tct^{\ol{\Phi}}$ and we will view $\tct^A$ as a (real) subbundle in
the complexification. In terms of this splitting the complex linear
extension of $\bJ$ is diagonalised. Concerning the tractor metric, we
first consider the unique Hermitian extension $\Cal H_{AB}$ of $\tilde
h_{AB}$ on $\tct_A$. Explicitly, $\Cal H_{AB}=\tilde h_{AB}-i\Bbb
J_{AB}$. This then extends to a complex bilinear form on
$\tct^A\otimes\Bbb C$. Writing $\tilde{h}$ and $\bJ$ for the complex
linear extensions of these tractors, in the matrix notation we have
$$
\bJ_B{}^A= \left(\begin{array}{cc} i\delta^\Phi_\Psi   & 0\\
0 & -i \delta^{\ol{\Phi}}_{\ol{\Psi}}
\end{array}\right) \quad \mbox{ and } \quad
 \tilde{h}_{AB}=\left(\begin{array}{cc}
    0                 & \Cal H_{\Phi\ol{\Psi}}\\
    \Cal H_{\ol{\Phi}\Psi} & 0
\end{array}\right)~,
$$
where $\Cal H_{\Phi\ol{\Psi}}$ is $\frac{1}{2}\Cal H_{AB}$, or more
accurately
$$
\frac{1}{2}\Cal H_{AB} = \left(\begin{array}{cc}
    0                 & \Cal H_{\Phi\ol{\Psi}}\\
    0 & 0
\end{array}\right),
$$
 and $\Cal H_{\ol{\Phi}\Psi}$ is the conjugate object. 

We write $\tZ_\Phi$ and $\tZ_\Phb$ for the holomorphic and
anti-holomorphic parts of the canonical section 
$X_A$ of the weighted 
conformal standard tractor bundle $\ct_A[1]$. That is 
$$
X_A= (\tZ_\Phi,\tZ_\Phb) .
$$
From this  and
$
\K_B=\bJ^A{}_B X_A.
$
it follows immediately that
$$
\K_B=(-i\tZ_\Phi,i\tZ_\Phb) \quad \Leftrightarrow \quad 
\K^B=(i\tZ^\Phi,-i\tZ^\Phb).
$$

For a choice of preferred scale on $\tM$, we write
$\frac{1}{2}\tY^\Phi$ and $\frac{1}{2}\tY^\Phb$ for the holomorphic
and anti-holomorphic parts of $Y^A$, i.e., $Y^A=
\frac{1}{2}(\tY^\Phi,\tY^\Phb).  $ (The normalisation on $\tY^\Phi$
means that we have $\tY^\Phi \tZ_\Phi=1$ which simplifies calculations
and is consistent with \cite{Gover-Graham}.)  It follows that
$$ 
\L^A= \tfrac{1}{ 2}(i\tY^\Phi,-i\tY^\Phb) \quad \Leftrightarrow
\quad \L_A= \tfrac{1}{2}(-i\tY_\Phi,i\tY_\Phb).
$$ 
Finally, we also have the complexified versions of the
$\tW$--projectors. The fact that $\tW$ is complex linear implies that
its complexification preserves the decomposition into holomorphic and
anti--holomorphic part, so we have
$\tW_A{}^a=(\tW_\Ph{}^\al,\tW_\Ph{}^\alb)$, and no combination of
barred and unbarred indices.

\newcommand{\tta}{{\tilde{\theta}}}
\newcommand{\tsi}{{\tilde{\sigma}}}
\newcommand{\Psb}{{\overline{\Psi}}}
\newcommand{\sib}{\overline{\sigma}}

\subsection{The case of a Fefferman space}\label{4.5}
If $\tilde M$ is the Fefferman space of a CR manifold $M$, then there
are several refinements of the above picture.  First observe that
there is a special class of (conformal) scales on $\tilde M$, namely
those coming from CR scales on $M$. A CR scale on $M$ simply is a
choice of positive contact form $\th$. As observed in \ref{2.3}, the
bundle $Q=TM/H$ naturally includes into $\ce(1,1)$. Now $\th$ defines
a linear map $Q\to\Bbb R$ and by complex linear extension a section of
the dual bundle $\ce(-1,-1)$. This section can be viewed as $U^{-1}$,
for a positive section $U$ of $\ce(1,1)$. Now as observed in \ref{3.5}
we may also view $U$ as a section of $\ce_{\Bbb C}[2]$, which is
easily seen to be in $\ce_+[2]$. Hence its square root can be used as
a scale $s$ on $\tilde M$.  We will also call such a conformal scale a
{\em CR scale}. By dint of context this should cause no confusion.  By
Theorem \ref{3.5}, sections of $\tce(1,1)$ are characterised among
sections of $\tce[2]$ by $\tilde{\Bbb D}_{\Bbb J}s=0$. So in fact CR
scales are exactly preferred scales in the sense introduced in the
last section and we can carry over the results that hold for these
scales. In particular in a CR scale $\k^a$ is a Killing field and
$Y_A\K^A=0=\L^AX_A$. In the subsequent calculations scales, when chosen, will 
be CR scales.

From the discussion of $\tilde{\Bbb D}$ in \ref{4.4} we see that
on density bundles $\tilde{\Bbb D}_{\Bbb J}$ differs from
$\tilde\nabla_\k$ by a multiple of $\tilde\nabla_a\k^a$, so in a
CR--scale the two operators coincide. By definition this implies that
in a CR scale, the double--$D$--operator $\tilde{\Bbb D}^\nabla_{\Bbb
J}$ on any weighted tractor bundle coincides with $\tilde\nabla_\k$.

\begin{prop*}
(1) The holomorphic and anti--holomorphic parts $\tilde Z^\Ph$ and $\tilde
Z^\Phb$ of $X^A\in\Ga(\tct^A_{\Bbb C}[-1])$ lie in the subspaces $\Ga(\Cal
T(1,0))$ respectively $\Ga(\Cal T^\Phb(0,1))$ and coincide with the
sections $Z^\Ph$ and $Z^\Phb$ from \ref{4.1}. 

\noindent
(2) For a choice $s$ of CR--scale, the sections $\tilde Y^\Ph$ and
    $\tilde Y^\Phb$ from \ref{4.4b} lie in the subspaces $\Ga(\Cal
    T(0,-1))$ and $\Ga(\Cal T(-1,0))$ and coincide with the sections
    $Y^\Ph$ and $Y^\Phb$ from \ref{4.2}
\end{prop*}
\begin{proof}
(1) By definition $2\tilde Z^\Ph=X^A-i\K^A$. Now $\k^a\tilde\nabla_a
    X^A=\k^a\Bbb Z^A{}_a=\K^A$. Since $\Bbb J$ is parallel, this shows
    that $\tilde Z^\Ph$ is an eigenvector for
    $\tilde\nabla_\k=\tilde{\Bbb D}^\nabla_{\Bbb J}$ with eigenvalue
    $i$, so by Proposition \ref{3.5} it lies in $\Ga(\Cal
    T(1,0))\subset\Ga(\tct\otimes\tce_{\Bbb C}[1])$. Viewed as an
    inclusion $\tce_{\Bbb C}[-1]\to \tct_{\Bbb C}$, $\tilde Z^\Ph$
    represents the complex linear extension of $X^A$. From \ref{2.5}
    we see that this represents the inclusion $\ce(-1,0)\to\Cal T$,
    which is given by $Z^\Ph$. The other statement follows in the same
    way. 

\noindent
(2) Next consider the tractor field $Y^A$ coming from a choice
$U\in\ce(1,1)\subset\tce_+[2]$ of CR scale on $\tM$. Applying the
fundamental derivative, we obtain $\tilde{\Bbb D}U\in
\tca^*\otimes\tce[2]$. Hence $U^{-1}\tilde{\Bbb
  D}U\in\tca^*\cong\tca$, where we use the trace--form for the last
identification. We claim that this is the \textit{grading element}
$E_s$ associated to the square root $s$ of $U$, i.e.~its eigenspaces
represent the splitting of $\tct$ according to $U$. Via the scale $U$,
the adjoint tractor bundle can be identified with $T\tilde
M\oplus\frak{co}(T\tilde M)\oplus T^*\tilde M$, and the sum of the
first and last part is orthogonal to the middle part with respect to
the trace--form. Also, the middle part splits into the orthogonal
direct sum of multiples of $\id_{T\tilde M}$ and $\frak{so}(T\tilde
M)$. In the scale determined by $U$, we have $\tilde\nabla U=0$ so
from the formula for $\tilde{\Bbb D}$ in \cite{luminy} (or
\cite{GoSrni99}) we get $\tilde{\Bbb D}_tU=0$ for $t\in T\tilde
M\oplus\frak{so}(T\tilde M)\oplus T^*\tilde M$. Hence
$U^{-1}\tilde{\Bbb D}U\in\Ga(\tca)$ is a multiple of $E_s$. By
definition $E_s\o E_s$ has trace two, so we can compute
$$
U^{-1}\tilde{\Bbb D}U=\tfrac{1}{2}E_sU^{-1}\tilde{\Bbb D}_{E_s}U=
\tfrac{1}{2}E_sU^{-1}2U=E_s,
$$ 
where we have used that the algebraic action of $E_s$ on $\ce[w]$ is
given by multiplication by $w$.  By Proposition \ref{3.6}, since $U$
lies in $\ce(1,1)\subset\tce[2]$, for a section $t$ of $\Cal A$, the
section $\tilde{\Bbb D}_tU$ lies in $\ce(1,1)$ and is equal to the CR
fundamental derivative $\Bbb D_tU$. Now one can play the same game as
above in the CR world to show that $U^{-1}\Bbb DU$, viewed as a
section of $\Cal A$, equals $\tfrac{1}{2}E_U$, where $E_U$ is the CR
grading element determined by the scale $U$. (The factor
$\tfrac{1}{2}$ is caused by the fact the $\pm 1$ eigenspaces of $E_U$
have each real dimension $2$, so applying the real trace--form to two
copies of $E_U$, one obtains $4$ rather than $2$.) This means that
$E_U$ is twice the component of $E_s$ in the decomposition of $\tca$
from \ref{3.4}, and hence twice the complex linear part of $E_s$.

Now consider the tractor field $Y^A$ determined by $s$. Viewed as a
projection $\tct\to\tct^1$, this is the projection onto the eigenspace
of $E_s$ with eigenvalue $1$. Explicitly,
$Y^A=\tfrac{1}{2}E_s\o(E_s+\id)$. A direct computation using
$E_s\o\Bbb J\o E_s=0$ shows that $\tfrac{1}{2}E_U\o(E_U+\id)$ is twice
the complex linear part of this projection. Hence decomposing $Y^A$
into holomorphic and anti--holomorphic parts, we obtain
$Y^A=\tfrac{1}{2}(Y^\Phi,Y^\Phb)$, for the weighted CR tractors
determined by the scale $U$ as in \ref{4.2}.
\end{proof}

\subsection{Relating the tangent bundles}\label{4.5a}
Our next task is to interpret the decomposition of the tangent bundle
from \ref{4.4a} in the special case of a Fefferman space $\pi:\tilde
M\to M$. The subbundle $V\subset T\tilde M$, spanned by $\k$, is the
vertical subbundle of $\pi$. 
\begin{lem*}
  Let $\pi:\tilde M\to M$ be a Fefferman space. Let $\th$ be a contact
  form on $M$ and consider the corresponding CR scale on $\tilde M$.
  Then we have $\k_a=2\pi^*\th$, $\tilde\nabla_a\k_b=2\pi^*d\th$ and the
  vector field $2\ell^a$ is the unique null lift of the Reeb field
  associated to $\th$.
  
  In particular, the subbundle $\tilde H\subset T\tilde M$ from
  \ref{4.4a} is exactly the preimage of the CR subbundle $H\subset
  TM$.
\end{lem*}
\begin{proof}
  The definition of $\tW$ in \eqref{tWdef} reads as 
$$
\bZ^A{}_{a}=\L^A k_a+ \K^A\ell_a+\tW^A{}_a ~.
$$
Since $\tilde{\nabla}_a X^A=\bZ^A{}_a$, we have $\k_a
=\K_A\bZ^A_{a}=\K_A\tilde{\nabla}_a X^A $. Since $X^A=(Z^\Phi,
Z^\Phb)$ and $K_A= (-iZ_\Phi,iZ_\Phb)$, this is
$-2iZ_\Phi\tilde{\nabla}_a Z^\Phi$.  Using $Z^\Phi Z_\Phi=0$, we see
that for any non-vanishing $\si\in \ce(-1,0)$ we have 
$Z_\Phi\tilde{\nabla}_a Z^\Phi = Z_\Phi\si^{-1}\tilde{\nabla}_a \si
Z^\Phi$. Now $\si Z^\Psi\in\Ga(\Cal T)\subset\Ga(\tct)$, so
Proposition \ref{3.2} implies that the one--form
$Z_\Phi\si^{-1}\tilde{\nabla}\si Z^\Ph$ is the pullback of
$Z_\Phi\si^{-1} \nabla\si Z^\Ph=Z_\Phi {\nabla}Z^\Ph$.  From the
formulae for the CR tractor connection we obtain $Z_\Phi\nd_a Z^\Phi =
i\theta$. Thus $k_a=2\pi^*\th_a$, so $\k^\perp$ coincides with the
preimage of $H$.

Next by part (1) of Proposition \ref{4.4}, $\tilde\nabla_a\k_b$ is
skew symmetric, so it coincides with $\tfrac{1}{2}$ times the exterior
derivative of $\k_a$, which equals $\pi^*d\th$. 

In view of these two results, the equations $\ell^a\k_a=1$ and
$\ell^b\tilde\nabla_b\k_a$ observed in Proposition \ref{4.4} imply
that the value of $2\ell^a$ at each point of $\tilde M$ projects onto
the Reeb field of $\tilde\th$. This pins down $\ell^a$ uniquely up to
adding $f\k^a$ for some smooth function $f$. But
$(\ell^a+f\k^a)(\ell_a+f\k_a)=2f$, which completes the proof. 
\end{proof}

Collecting the results, we see that for a Fefferman space $\pi:\tilde
M\to M$, the filtration $V\subset\k^\perp\subset T\tilde M$ from
\ref{4.4a} has the form $\ker(T\pi)\subset T\pi^{-1}(H)\subset T\tilde
M$. Note further that the resulting identification of $H$ with $\tilde
H/V$ is compatible with the complex structures on both bundles, since
they were both induced from the complex structure on the tractor
bundle.

Next given a choice of CR scale, we can explicitly identify sections of
the CR subbundle $H\to M$ with sections of the corresponding subbundle
$\tilde H=\k^\perp\cap\ell^\perp\subset T\tilde M$. Since it will
be  useful later, we do this in a complexified picture and in a
weighted version.

\begin{prop*}
  Let $\pi:\tilde M\to M$ be a Fefferman space, and fix some CR scale.
  Let $\tilde H\otimes\Bbb C=\tilde H^\al\oplus\tilde H^\alb$ be the
  decomposition of the complexification of $\tilde
  H=\k^\perp\cap\ell^\perp$ into holomorphic and anti--holomorphic
  parts. Then for arbitrary weights $w$ and $w'$, sections of
  $H^\al(w,w')$ are in bijective correspondence with sections $\xi$ of
  $\tilde H^\al[w+w']$ such that $\tilde\nabla_\k\xi=(w-w'+1)i\xi$.
  Likewise, sections of $H^\alb(w,w')$ are in bijective correspondence
  with sections $\xi$ of $\tilde H^\alb[w+w']$ such that
  $\tilde\nabla_\k\xi=(w-w'-1)i\xi$.
\end{prop*}
\begin{proof}
  Let us first treat the real subbundles $H$ and $\tilde H$. The flow
  lines of $\k^a$ are exactly the fibres of $\pi$. From this one
  easily concludes that a vector field $\xi\in\frak X(\tM)$ is
  projectable if and only if the Lie derivative $\Cal L_\k\xi$ is
  vertical, and thus $\Cal L_\k\xi=\ell(\Cal L_\k\xi)\k$ (where we
  view $\ell$ as a 1-form). If $\xi$ is a section of the subbundle
  $\k^\perp\cap\ell^\perp$, then $\ell(\xi)=0$, and hence $\ell(\Cal
  L_\k\xi)=-(\Cal L_\k\ell)(\xi)$.  In \ref{4.4} we have observed that
  $\ell_a=\tilde\Rho_{ab}\k^b$. Also there we noted that since, in a CR scale, 
$\k$ is Killing, 
we have $\Cal L_\k\ell=0$ and it follows that a section $\xi$ of
  $\k^\perp\cap\ell^\perp$ is projectable if and only if $\Cal
  L_\k\xi=0$. Since $\tilde\nabla$ is torsion free, we get $\Cal
  L_\k\xi=\tilde\nabla_\k\xi-\tilde\nabla_\xi\k$. Thus, sections of
  $H$ are in bijective correspondence with sections $\xi$ of $\tilde
  H$ such that $\tilde\nabla_\k\xi=\tilde\nabla_\xi\k$.
  
  Now, from Proposition \ref{4.4}, $\tilde\nabla_\xi\k$ is the complex
  structure on $\tilde H$ applied to $\xi$, so on sections of $\tilde
  H^\al$ this coincides with $i\xi$, and on sections of $\tilde H^\alb$
  it coincides  with $-i\xi$. Thus we obtain the result for $w=w'=0$.
 
  To conclude the proof, we observe that a CR scale is a preferred
  scale and hence is killed by $\tilde{\Bbb D}_{\Bbb J}$. Since powers
  of this scale are used to trivialise density bundles, we conclude
  that $\tilde{\Bbb D}_{\Bbb J}=\tilde\nabla_\k$ on density bundles.
  Using this, the general result immediately follows from Proposition
  \ref{3.5}.
\end{proof}

\subsection{Computing a metric from the conformal class}\label{4.8}
A choice of CR scale determines a metric from the canonical conformal
class on a Fefferman space. We are now ready to compute this explicitly. 
 Choose
a contact form $\th$ on $M$, and let $U\in\ce(1,1)$ be the
corresponding CR scale. Choose a section $\si$ of
$\ce(-1,0)\subset\tce_{\Bbb C}[-1]$ such that
$\si\overline{\si}=U^{-1}$. Locally, we define a smooth function
$\ga_{\si}:\tilde M\to\Bbb R$ by requiring that $\tilde\si=\sigma
e^{i\gamma_\si}$ is a positive real section of $\tce[-1]$. By
definition, this implies that $\tilde\si=s^{-1}$, where $s$ is the CR
scale determined by $U$. Using that the points in a Tiber of
$\pi:\tilde M\to M$ determine how the real line $\tct^1$ sits inside
the complex line $\Cal T^1$, one easily verifies that $\ga_\si$ defines
a local coordinate for each fibre.

On the other hand, recall that the real part of the Levi form defines
a non--degenerate bundle metric on the CR subbundle $H$. We can
uniquely extend this to a (degenerate) bundle metric $L:TM\x TM\to\Bbb
R$ by requiring that the Reeb vector field inserts trivially into $L$.
This is called the \textit{degenerate Levi metric}.

\begin{prop*} 
  Let $\pi:\tilde M\to M$ be a Fefferman space, $\th$ a contact form,
  $U\in\Ga(\ce(1,1))$ the corresponding CR scale. Choose
  $\si\in\Ga(\ce(1,0))$ such that $\si\bar\si=U$ and consider the
  (local) one--form $\ga_\si$ defined above. 
  
  Then the one--form $\tau:=\ell_a$
$$
\tau=-\frac{i}{2}\pi^*(\si^{-1}\nd \si- \sib^{-1}\nd \sib) -\frac{1}{n+2}
\tta \pi^*(P) + d \gamma_\si ~
$$
depends only on $\th$ and the metric $g_\th$ in the conformal class
corresponding to the CR scale $s$ determined by $U$ is given by
$$
g_\theta=  \pi^*L+4\tau \odot \pi^*\th
$$
where $L$ is the degenerate Levi metric.
\end{prop*}
\begin{proof}
  We start by computing $\ell_a=\L_B\bZ^B{}_a$. Since
  $\bZ^B{}_a=\tilde{\nabla}_a X^B$, we have to calculate
  $\L_B\tilde{\nabla}_a X^B$. Recall that
$$
\tilde{\nabla}_a X^B=( \tilde{\nabla}_a Z^\Psi, \tilde{\nabla}_a
Z^\Psb ) \quad \mbox{and} \quad \L_B= \frac{1}{2}(-iY_\Psi,iY_\Psb). 
$$ 

\noindent
The section $\tsi=s^{-1}$ is parallel for $\tilde\nabla$, which
implies $\si\tilde{\nabla}_aZ^\Psi=e^{-i\gamma_\si}\tilde{\nabla}_a
\tsi Z^\Psi$. We apply the Leibniz rule to rewrite this and obtain
$(\tilde{\nabla}_a \si Z^\Psi + i\si Z^\Psi \tilde{\nabla}_a
\gamma_\si)$.  So $ \tilde{\nabla}_a Z^\Psi = \si^{-1}\tilde{\nabla}_a
\si Z^\Psi + i Z^\Psi \tilde{\nabla}_a \gamma_\si $. Contracting this
with $iY_\Psi$ we get $i\si^{-1}Y_\Psi\tilde{\nabla}_a \si
Z^\Psi-\tilde{\nabla}_a\gamma_\si $. Now $\si Z^\Psi\in\Ga(\Cal
T)\subset\Ga(\tct)$, which implies that the one--form
$i\si^{-1}Y_\Psi\tilde{\nabla} \si Z^\Psi$ is the pullback of
$i\si^{-1}Y_\Psi\nd \si Z^\Psi$. But by the Leibniz rule and the
formulae for $\nd Z^\Phi$ we have
 $$
i\si^{-1}Y_\Psi\nd \si Z^\Psi =  i\si^{-1}\nd \si +\tfrac{1}{n+2} \theta P
$$ 
Thus $ iY_\Psi\tilde{\nabla}Z^\Psi=\pi^*(i\si^{-1}\nd \si
+\frac{1}{n+2} \theta P) -\nd \gamma_\si $ and averaging this with
its conjugate brings us to
$$
\tau = -\tfrac{i}{2}\pi^*(\si^{-1}\nd \si- \sib^{-1}\nd \sib) 
-\tfrac{1}{n+2} \tta \pi^*(P) + d \gamma_\si  ,
$$ 
A simple computation
shows that $\tfrac{i}{2}\pi^*(\si^{-1}\nd \si- \sib^{-1}\nd \sib)- d
\gamma_\si$ depends only on $\th$ and not on the choice of $\si$.

In view of the formula \nn{metric} for $g_{ab}$ it remains to discuss
the quantity $\h_{ab}$ occurring there. Since it annihilates $\k^a$,
at a point $x\in\tilde M$, it descends to $T_{\pi(x)}M$. Since
$\ell^a\h_{ab}=0$, this descended quantity annihilates the Reeb field.
On the other hand, from the construction of $\h_{ab}$ via the tractor
metric it follows that the restriction to $H_{\pi(x)}$ coincides with
the real part of the Levi form. But this immediately implies that
$\h_{ab}$ is the pullback of $L$.
\end{proof}

\newcommand{\tL}{{\tilde L}}

\subsection{Relating preferred connections}\label{4.8a}
Let $\pi:\tilde M\to M$ be a Fefferman space. Choosing a contact form
$\th$ (on $M$) we get the Webster--Tanaka connection on $M$ and an induced CR
scale on $\tilde M$. We next want to compare the Levi--Civita
connection associated to the latter with the downstairs Webster--Tanaka connection.

\begin{prop*}
  Let $\pi:\tilde M\to M$ be a Fefferman space, $\th$ a choice of
  contact form on $M$ and $s$ the corresponding CR scale on $\tilde
  M$. Let us denote by $\nabla$ the corresponding Webster--Tanaka
  connections (if necessary coupled to CR tractor connections) and by
  $\tilde\nabla$ the corresponding Levi--Civita connections (if necessary coupled to
  conformal tractor connections). By $\nabla^H$ we denote the
  restriction of $\nabla$ to directions in $H\subset TM$.

\noindent
(1) For a complex conformal weighted tractor bundle $\tcw[w+w']$, and
any section $f\in\Ga(\Cal W(w,w'))\subset\Ga(\tcw[w+w'])$, $\Bbb
I_a^c\tilde\nabla_cf$ descends to $\nabla^Hf$, and
$2\ell^a\tilde\nabla_a f$ descends to
$\nabla_0f+\tfrac{i(w-w')}{n+2}\Rho f$.

\noindent
(2) For $\xi^a\in\Ga(H^\al(w,w'))\subset\Ga(\tilde H^a[w+w'])$, $\Bbb
I^a_c\Bbb I_b^d\tilde\nabla_d\xi^c$ descends to $\nabla^H\xi^\al$, and
$2\Bbb I^a_c\ell^b\tilde\nabla_b\xi^c$ descends to
$\nabla_0\xi^\al+i\Rho^\al{}_\be\xi^\be+\tfrac{i(w-w')}{n+2}\Rho
\xi^\al$.

\noindent
(3) Consider the extension of the projector/injector $\tW_A{}^a$
associated to $s$ to the complexified tractor bundle. Then the
decomposition into holomorphic and anti--holomorphic part has the form
$\tW_A{}^a=(\tW_\Ph{}^\al,\tW_\Phb{}^\alb)$, and the two components
descend to the CR--objects $W_\Ph{}^\al\in\Ga(\ce_\Ph{}^\al(-1,0))$
and $W_\Phb{}^\alb\in\Ga(\ce_\Phb{}^\alb(0,-1))$ from \ref{4.2}.
\end{prop*}
\begin{proof}
  (1) For unweighted tractor bundles, $\tilde\nabla$ and $\nabla$ are
  just tractor connections, so the results follow from Proposition
  \ref{3.2}. By the Leibniz rule it therefore suffices to prove the
  result for densities. Assume first, that $f\in\Ga(\ce(-1,0))$.
  Working on $M$, we can use $Y_\Ph Z^\Ph=1$ to compute
$$
\nabla f=\nabla Y_\Ph f Z^\Ph=Y_\Ph\nabla fZ^\Ph+fZ^\Ph\nabla Y_\Ph. 
$$
From the formulae for the components of $\nabla Y_\Ph$ in \ref{4.3}
we see that $Z^\Ph\nabla^HY_\Ph=0$ and
$Z^\Ph\nabla_0Y_\Ph=\tfrac{i}{n+2}\Rho$. On the other hand, $f Z^\Ph$
is an unweighted tractor, so we know that $\tilde\nabla fZ^\Ph$
descends to $\nabla fZ^\Ph$. Using the Leibniz rule once more, we get
$$
Y_\Ph\Bbb I_b^d\tilde\nabla_d fZ^\Ph=\Bbb I_b^d\tilde\nabla_d
f+Y_\Ph f\Bbb I_b^d\tilde\nabla_d Z^\Ph,
$$
and this descends to $\nabla^Hf$. By \eqref{conf-connection}, $\Bbb
I_b^d\tilde\nabla_dX^A=\Bbb I_b^d\Bbb Z^A{}_d=\tW^A{}_b$. Since
$Z^\Ph$ lies in the complex subspace generated by $X^A$, we see that
$\Bbb I_b^d\tilde\nabla_d Z^\Ph$ lies in the complex subspace
generated by $\tW^A{}_b$. Since both $Y_A$ and $L_A$ hook trivially
into $\tW^A{}_b$ (and hence into any element of that complex
subspace), we conclude that $\Bbb I_b^d\tilde\nabla_d f$ descends to
$\nabla^H f$. Passing to powers of $f$ and $\overline{f}$, we see that
this holds for arbitrary densities.

To deal with $\nabla_0$, recall from Lemma \ref{4.5a} that $2\ell^a$
is a lift of the Reeb field. Since $fZ^\Ph$ is an unweighted tractor,
$2Y^\Ph\ell^a\tilde\nabla_a f Z^\Ph$ descends to
$$
Y^\Ph\nabla_0 fZ^\Ph=\nabla_0f-\tfrac{i}{n+2}\Rho f. 
$$ 
Since $\tilde\nabla_a X^A=\Bbb Z^A{}_a$, we see that
$\ell^a\tilde\nabla_a Z^\Ph$ lies in the complex subspace
generated by $\ell^a\Bbb Z^A{}_a=\L^A$, which implies that $\tilde
Y_\Ph\ell^a\tilde\nabla_a Z^\Ph=0$, so $2\ell^a\tilde\nabla_a
f$ descends to $\nabla_0f-\tfrac{i}{n+2}\Rho f$. Using powers of $f$
and its conjugate we obtain the general formula.

\noindent
(3) In the proof of Proposition \ref{4.5a} we have noted that
$\tilde{\Bbb D}_{\Bbb J}=\tilde\nabla_\k$ on density bundles. By
definition, this implies $\tilde{\Bbb D}^\nabla_{\Bbb
  J}=\tilde\nabla_\k$ on weighted tractor bundles.

We first claim that $\k^b\tilde\nabla_b\tilde{W}_A{}^a=0$. By
\eqref{conf-connection}, $\k^b\tilde\nabla_b\Bbb
Z_A{}^a=-\ell^aX_A-\k^aY_A$. On the other hand, using that
$L_A=Y_B\Bbb J^B{}_A$ and $\k^b\tilde\nabla_b\k^a=0$ we get
$\k^b\tilde\nabla_b\L_A\k^a=-Y_A\k^a$. Likewise,
$\k^b\tilde\nabla_b\K_A=-X_A$, so inserting the definition (\ref{tWdef}) of $\tilde{W}_A{}^a$
we obtain $\k^b\tilde\nabla_b\tilde{W}_A{}^a=-\K_A\k^b\tilde\nabla_b\ell^a$,
which vanishes by part (4) of Proposition \ref{4.4}.

For $s^A\in\Ga(\tct_{\Bbb C}[1])$ we therefore get
$\k^c\tilde\nabla_cs^A\tW_A{}^a=\tW_A{}^a\k^c\tilde\nabla_cs^A$.  If
$s^A$ lies in the subspace $\Ga(\Cal T^\Ph(1,0))$, then
$\k^c\tilde\nabla_cs^A=is^A$ and hence
$\k^c\tilde\nabla_cs^A\tW_A{}^a=is^A\tW_A{}^a$. By Proposition
\ref{4.8} we see that $s^A\tW_A{}^a\in\Ga(H^\al)\subset\Ga(\tilde
H^\al)$. Likewise, if $s^A\in \Ga(\Cal T^\Phb(0,1))$ we get
$s^A\tW_A{}^a\in\Ga(H^\alb)$, so we see that the holomorphic and
anti--holomorphic parts of $\tW_A{}^a$ descend as required. Since the
appropriate parts of the complementary projections $X_A$ and $Y_A$
descend to their CR analogs, these descended sections must coincide
with $W_\Ph{}^\al$ and $W_\Phb{}^\alb$.

\noindent
(2) For a section $\xi^a$ of (the complexification of) $\tilde H$ we
    compute
    $$
    \Bbb I^a_c\tilde\nabla_b
    \xi^c=\tW_A{}^a\tW^A{}_c\tilde\nabla_b\xi^c=\tW_A{}^a\tilde\nabla_b
    \tW^A{}_c\xi^c-\tW_A{}^a\xi^c\tilde\nabla_b\tW^A{}_c.
$$
By \eqref{tWdef}, $\tilde\nabla\tW^A{}_c-\tilde\nabla\Bbb Z^A{}_c$
is a sum of terms which contain one of the four elements $\K^A$,
$\L^A$, $\k_c$, or $\ell_c$ undifferentiated. But the first two are
killed by contraction with $\tW_A{}^a$, while the last two are annihilated  by contraction into $\xi^c$.
Finally, $\tW_A{}^a\tilde\nabla_b\Bbb Z^A{}_c=0$ by
\eqref{conf-connection}, so we obtain
$$
\Bbb I^a_c\tilde\nabla_b \xi^c=\tW_A{}^a\tilde\nabla_b \tW^A{}_c\xi^c.
$$
If $\xi^a$  lies in $\Ga(H^\al)\subset\Ga(\tilde
H\otimes\Bbb C)$, then by part (3), $\tW^A{}_c\xi^c$ lies in $\Ga(\Cal
T_\Ph(-1,0))\subset\Ga(\tct_{\Bbb C}[-1])$, so we can apply part (1).
Together with part (3) we see that $\Bbb I^a_c\Bbb
I_b^d\tilde\nabla_d\xi^c$ descends to $W_\Ph{}^\al\nabla^H
W^\Ph{}_\be\xi^\be$.  Applying the Leibniz rule, and using that the
formulae in \ref{4.3} show that $W_\Ph{}^\al\nabla^H W^\Ph{}_\be=0$, we
conclude that $\Bbb I^a_c\Bbb I_b^d\tilde\nabla_d\xi^c$ descends to
$\nabla^H\xi^\al$.

Still assuming that $\xi^a\in \Ga(H^\al)\subset\Ga(\tilde H\otimes\Bbb
C)$, we see that $2\Bbb I^a_c\ell^b\tilde\nabla_b \xi^c$ descends to
$W_\Ph{}^\al\nabla_0 W^\Ph{}_\be\xi^\be+\tfrac{i}{n+2}\Rho\xi^\al$.
The formulae for $\nabla W$ in \ref{4.3} show that
$$
W_\Ph{}^\al\tilde\nabla_0
W^\Ph{}_\be=i\Rho_{\be}{}^{\al}-\tfrac{i}{n+2}\Rho\delta^\al_\be,
$$
so $2\Bbb I^a_c\ell^b\tilde\nabla_b \xi^c$ descends to
$\nabla_0\xi^\al+i\Rho_{\be}{}^{\al}\xi^\be$. Since we have
established the formulae for densities already in part (1), this
completes the proof.
\end{proof}

\subsection{Conformal Killing fields on Fefferman spaces}\label{4.8b}
We can now make the decomposition of conformal Killing fields on a
Fefferman space from Theorem \ref{3.4} explicit.
\begin{thm*}
  Let $M$ be a CR manifold with Fefferman space $\tilde M$, let
  $v^a$ be a conformal Killing field on $\tilde M$, and fix a choice
  of preferred scale.

\noindent
(1) $v^b\tilde\nabla_b\k^a-\k^b\tilde\nabla_b v^a$ is a conformal
    Killing field on $\tilde M$ which inserts trivially into the
    tractor curvature $\Om_{ab}{}^C{}_D$. 

\noindent
(2) The vector field
$u^a:=v^a-(\k^c\tilde\nabla_cv^d)\tilde\nabla_d\k^a+\k_bv^b\ell^a$ on
$\tilde M$ descends to an infinitesimal CR automorphism of $M$.
Further $\k_bv^b$ descends to a smooth function on $M$ from which this
infinitesimal automorphism can be recovered by a CR--invariant
differential operator.

\noindent
(3) Define $w^a:=\Bbb
    I^a_cv^c+(\k^b\tilde\nabla_bv^d)(\tilde\nabla_d\k^a)$. Then the
    section $w^a-iw^c\tilde\nabla_c\k^a$ of $\tilde H\otimes\Bbb C$
    descends to a section $w^\al$ of $\ce^\al(-1,1)$ which satisfies
    $\nabla_\al w^\be=\delta_\al^\be\nabla_\ga w^\ga$ as well as
    $\nabla^\al w^\be=-\nabla^\be w^\al$.  
\end{thm*}
\begin{proof}
As in  the proof of Proposition \ref{4.4} we put
$$
V^B:=\Bbb Z^B{}_av^a-\tfrac{1}{2n+2}X^B\tilde\nabla_a v^a
$$ 
and consider the adjoint tractor field $s_A{}^B=\tfrac{1}{2n+2}D_AV^B$
associated to $v^a$.

\noindent
(1) By Lemma \ref{3.4}, $\{s,\Bbb J\}$ is a parallel section of
    $\tca$, so the underlying vector field is conformal Killing and
    inserts trivially into the Cartan curvature. We can compute this
    underlying vector field as 
$$
X^A\Bbb Z_B{}^a(s_A{}^C\Bbb J_C{}^B-\Bbb J_A{}^C s_C{}^B). 
$$
Using part (5) of Proposition \ref{4.4} and formula \eqref{tractor-D},
this expands as 
$$
V^C(Y_C\k^a+\Bbb Z_C{}^c\tilde\nabla_c\k^a+X_C\ell^a)-\Bbb
Z_B{}^a\k^c\tilde\nabla_cV^B, 
$$
which easily leads to the required expression. 

\noindent
(2) By Theorem \ref{3.4}, we can form the complex linear part of
$s_A{}^B$ and add an appropriate multiple of $\Bbb J_A{}^B$ to obtain
an element of $\Ga(\Cal A)\subset\Ga(\tca)$ which defines an
infinitesimal CR automorphism of $M$. In particular, the corresponding
conformal Killing field on $\tilde M$ descends to that infinitesimal
CR automorphism on $M$. Since the multiple of $\Bbb J_A{}^B$ just
contributes a multiple of $\k^a$ to the underlying vector field, we
can ignore it in the computation. The conformal Killing field
underlying the complex linear part of $s_A{}^B$ can be computed as 
$$
\tfrac{1}{2} X^A\Bbb Z_B{}^a(s_A{}^B-\Bbb J_A{}^Cs_C{}^D\Bbb J_D{}^B).
$$
Using formula \eqref{tractor-D} and part (5) of Proposition \ref{4.4}
this is easily evaluated directly, and one obtains
$$
\tfrac{1}{2}\left(u^a+(\ell_bv^b+
\tfrac{1}{2n+2}\k^c\tilde\nabla_c\tilde\nabla_b
v^b)\k^a\right).
$$
Hence we see that $u^a$ descends to an infinitesimal CR
automorphism on $M$. From the definition of $u$ we immediately see
that $\k_au^a=2\k_a v^a$. In view of Lemma \ref{4.5a}, we obtain the
section of $TM/H$ induced by this infinitesimal automorphism by
multiplying the Reeb field by the function $\tfrac{1}{2}\k_av^a$. Now
an infinitesimal CR automorphism on $M$ can be recovered by applying
an invariant differential operator to its projection to $TM/H$, see
\cite[3.4]{deformations}.

We can also verify the last two facts directly: Using Proposition
\ref{4.4} we get
$\k^b\tilde\nabla_b\k_cv^c=\k^b\k^c\tilde\nabla_bv_c$, and this
vanishes since the symmetric part of $\tilde\nabla_bv_c$ is pure trace
and $\k$ is isotropic. Thus, the function $\k_cv^c$ descends to $M$.
To recover the infinitesimal CR automorphism from this function, it
suffices to recover $\Bbb I_b^a u^b=\Bbb I_b^av^b-(\k^c\tilde\nabla_c
v^d)\tilde\nabla_d\k^a$. To do this, we use that $v$ is conformal
Killing and $\k$ is Killing to compute
$$
\tilde\nabla^c\k_bv^b=-\k^b\tilde\nabla_bv^c+
\tfrac{1}{n+1}\k^c\tilde\nabla_bv^b-v^b\tilde\nabla_b\k^c,
$$
which implies that $\Bbb I_b^a u^b=(\tilde\nabla^c\k_b
v^b)(\tilde\nabla_c\k^a)$, since $\tilde{\nabla}_b\k^c$ is the
complex structure on $\tilde{H}$ and hence
$(\tilde\nabla_b\k^c)(\tilde\nabla_c\k^d)=-\Bbb I_b^d$.

\noindent
(3) As before let $s_A{}^B$ be the adjoint tractor field corresponding
to the conformal Killing field $v^a$. From Theorem \ref{3.4} we know
that the conjugate linear part of $s_A{}^B$ descends to a section of
$\Cal E^{[\Ph\Ps]}$, and that there is a canonical projection from
that bundle to its irreducible quotient $\Cal E^\al(-1,1)$. Applying
this projection to the conjugate linear part, we obtain a section in
the kernel of two CR invariant operators, and from the Remark after
Theorem \ref{3.4} we know that this corresponds to the two claimed
properties of $w^\al$. From the definitions, one easily concludes that
the projection $\ce^{[\Ph\Ps]}\to\ce^\al(-1,1)$ is explicitly given by
$Z_\Ph W_\Ps{}^\al$. Now a direct computation shows that
$$
w^a=X^A\tW_B{}^a(s_A{}^B+\Bbb J_A{}^C s_C{}^D\Bbb J_D{}^B), 
$$
so the section of $\ce^\al(-1,1)$ in question can be computed as the 
holomorphic part $w^a-iw^c\tilde\nabla_c\k^a$ of $w^a$.

Alternatively, one may verify all the claims by direct computations
along the following lines: We can write
$$
w^a=(\k^b\tilde\nabla_b v^d-v^b\tilde\nabla_b \k^d)(\tilde\nabla_d\k^a)
$$
and differentiate this formula. To expand these derivatives, one has
to use the differential consequences of the conformal Killing equation
as detailed in \cite{Gopowers}: Putting
$\rho_a:=-\tfrac{1}{2n+2}\nabla_a\nabla_cv^c-\tilde\Rho_{ac}v^c$, one
has
$$
  \tilde\nabla_a\tilde\nabla_b v^c=g_{ab}\rho^c-\delta_a^c\rho_b-
\tilde\Rho_{ab}v^c+\tilde\Rho_a{}^cv_b
+C_b{}^{cd}{}_av_d-\rho_a\delta_b^c-\tilde\Rho_{ad}v^d\delta_b^c,
$$
where $C_{ab}{}^c{}_d$ is the Weyl curvature. Moreover, in the
corresponding equation for $\k$ instead of $v$, one may replace
$\rho_a$ by $-\ell_a$ and the last three summands vanish. 

Using these identities, one easily verifies directly that
$\k^b\tilde\nabla_bw^a=-w^b\tilde\nabla_b\k^a$. This immediately
implies that the section $w^a-iw^c\tilde\nabla_c\k^a$ of $\tilde
H\otimes\Bbb C$ is an eigenvector for the operator
$\k^b\tilde\nabla_b$ with eigenvalue $-i$. By Proposition \ref{4.5a},
this implies that it descends to a section $w^\al$ of $\ce^\al(-1,1)$.

Next, one computes $\Bbb I^a_c\Bbb I^d_b\tilde\nabla_d w^c$ by first
expanding the expression for $\nabla_dw^c$ and ignoring those terms
which have a free index on either $\k$ or $\ell$. This leads to
\begin{align*}
\Bbb I^a_c\Bbb I^d_b\tilde\nabla_d w^c=&\Bbb I^a_c\Bbb
I^d_b\tilde\nabla_d v^c+(\tilde\nabla_b\k^c)(\tilde\nabla_c v^d)
(\tilde\nabla_d\k^a)\\
+&\Bbb I^a_b\k^c\ell_d\tilde\nabla_cv^d-(\k^c\rho_c+\ell_c
v^c)\tilde\nabla_b\k^a.  
\end{align*}
The first and second line of this formula exactly are the conjugate
linear part respectively the complex linear part of the resulting
endomorphisms of $\tilde H$. The second line is evidently a complex
multiple of the identity. On the other hand, the fact that the
symmetrisation of $\tilde\nabla_c v^d$ is pure trace easily implies
that the first line is skew symmetric. Using part (2) of Proposition
\ref{4.8a}, this easily implies the result.
\end{proof}

\subsection{Relating double--$D$'s}\label{4.9}
Consider the conformal double--$D$--operator $\tilde{\Bbb D}^{\nabla}$
as an operator mapping sections of a weighted complex conformal
tractor bundle $\tcw[w+w']$ to sections of $\tca_{\Bbb
  C}\otimes\tcw[w+w']$, see \ref{3.7}. Hence we denote the operator by
$(\tilde{\Bbb D}^\nabla)_A{}^B$. Having a complex conformally natural
bundle $\tcw$, we can use the splitting of $\tca_{\Bbb C}$ from
\ref{3.7} to obtain
\begin{equation}\label{fund-split}
(\tilde{\Bbb D}^{\nabla})_A{}^B=\begin{pmatrix}(\tilde{\Bbb
D}^{\nabla})_{\Ph}{}^{\Ps} & (\tilde{\Bbb
D}^{\nabla})_{\Phb}{}^\Ps\\ 
(\tilde{\Bbb D}^{\nabla})_{\Ph}{}^{\Psb}& (\tilde{\Bbb
D}^{\nabla})^{\Phb}{}_{\Psb}\end{pmatrix}
\end{equation}
The fact that $(\tilde{\Bbb D}^\nabla)_A{}^B$ is skew symmetric
implies that the component in the bottom left corner is the negative
transpose of the one in the top right corner, so the notation is
consistent.

From Theorem \ref{3.7} we know that $(\tilde{\Bbb
D}^\nabla)_{\Ph}{}^{\Ps}$ restricts to the CR--double--$D$--operator
$$
(\Bbb D^{\nabla})_{\Ph}{}^{\Ps}:\Ga(\Cal W(w,w'))\to\Ga(\Cal
A_{\Bbb C}\otimes \Cal W(w,w')).
$$
On the other hand, the off diagonal components descend to operators 
\begin{gather*}
(\Bbb D^{\nabla})_{\Ph\Ps}:\Ga(\Cal W(w,w'))\to\Ga(\Cal
  E_{[\Ph\Ps]}\otimes \Cal W(w-1,w'+1))\\
(\Bbb D^{\nabla})_{\Phb\Psb}:\Ga(\Cal W(w,w'))\to\Ga(\Cal
  E_{[\Phb\Psb]}\otimes \Cal W(w+1,w'-1)).
\end{gather*}

\begin{thm*}
Consider a weighted complex CR tractor bundle $\Cal W(w,w')$. Then the
operators $\Bbb D$ constructed above are explicitly given by
\begin{align*}
\Bbb D_{\Ph\Ps}f=&2wZ_{[\Ps}Y_{\Ph]}f+
2Z_{[\Ps}W_{\Ph]}{}^\al\nabla_\al f\\ 
\Bbb D_{\Phb\Psb}f=&2w'Z_{[\Psb}Y_{\Phb]}f+
2Z_{[\Psb}W_{\Phb]}{}^\alb\nabla_\alb f\\
\Bbb D_{\Ph\Psb}f=&w Z_\Psb Y_\Ph f-w'Z_\Ph Y_\Psb f+Z_\Psb
W_\Ph{}^\al\nabla_\al f-Z_\Ph W_\Psb{}^\alb\nabla_\alb f\\
-&Z_\Ph Z_\Psb(i\nabla_0f+\tfrac{w'-w}{n+2}\Rho f).
\end{align*}
\end{thm*}
\begin{proof}
  By \cite[section 4]{GoSrni99}, the conformal
  double $D$--operator on $\tcw[w+w']$ is given by
$$
  \tilde{\Bbb D}^{\nabla}_{AB}f=2(w+w')X_{[B}Y_{A]}f+2X_{[B}\Bbb
  Z_{A]}{}^a\tilde\nabla_a f,
  $$
  see also \cite{Goadv}.  Now let us insert $\Bbb
  Z_A{}^a=\tW_A{}^a+\L_A\k^a+\K_A\ell^a$. From the proof of part (3)
  of Proposition \ref{4.8a}, we know that, in a CR scale,
  $\tilde\nabla_\k$ coincides with $\tilde{\Bbb D}^\nabla_{\Bbb J}$ on
  weighted tractor bundles, so $\k^a\tilde\nabla_a f=i(w-w')f$. Using
  this, we can rewrite $\frac{1}{2}\tilde{\Bbb D}^\nabla_{AB}f$ as
$$ (w+w')X_{[B}Y_{A]}f+X_{[B}\tW_{A]}{}^a\tilde\nabla_a
f+i(w-w')X_{[B}\L_{A]}f+X_{[B}\K_{A]}\ell^a\tilde\nabla_a f.
$$
For each term occurring in this decomposition, we understand explicitly the
decomposition into holomorphic and anti--holomorphic parts. In
particular, the components of $(w+w')Y_A+i(w-w')L_A$ descend to
$(wY_\Ph,w'Y_\Phb)$, compare with \ref{4.4b}. Further, the components
of $X_B$, $\K_A$, and $\tW_A{}^a$ descend to $(Z_\Ps,Z_\Psb)$,
$(-iZ_\Ph,iZ_\Phb)$, and $(W_\Ph{}^\al,W_\Phb{}^\alb)$, respectively.
Inserting this, decomposing and using Proposition \ref{4.8a}, the
claimed formulae follow.
\end{proof}

\subsection{Relating tractor--$D$'s}\label{4.10}
To complete our picture, it remains to interpret the conformal
Rho--tensor $\tilde\Rho$ in terms of CR--data.
\begin{prop*}
  Let $\pi:\tilde M\to M$ be a Fefferman space, and consider the CR
  scale on $\tilde M$ induced by a choice of contact form on $M$. Then
  the complex bilinear extension of the conformal Rho--tensor is given
  by
  \begin{align*}
    \tilde\Rho_{ab}=&\ell_a\ell_b-\tfrac{1}{4}S\k_a\k_b+i\k_{(a}\Bbb
      I_{b)}^\al T_\al-i\k_{(a}\Bbb I_{b)}^\alb T_\alb\\
+&\tfrac{i}{2}\Bbb I_a^\al\Bbb I_b^\be A_{\al\be}-\tfrac{i}{2}\Bbb
      I_a^\alb\Bbb I_b^\beb A_{\alb\beb}+\tfrac{1}{2}\Bbb I_a^\beb\Bbb
      I_b^\al\Rho_{\al\beb}+\tfrac{1}{2}\Bbb I_a^\be\Bbb
      I_b^\alb\Rho_{\alb\be}.
  \end{align*}
In particular, $\tilde\Rho_a{}^a=\Rho_\al{}^\al$. 
\end{prop*}
\begin{proof}
  By \eqref{conf-connection}, we have $\tilde\Rho_{ab}=\Bbb
  Z^A{}_b\tilde\nabla_aY_A$. Now we can decompose
$$
\tilde\nabla_aY_A=(\ell_a\k^c+\k_a\ell^c+\Bbb I_a^c)\tilde\nabla_c Y_A.
$$
By Proposition \ref{4.5}, $Y_A=\tfrac{1}{2}(Y_\Ph,Y_\Phb)$ and the
components lie in $\Ga(\Cal T_\Ph(-1,0))$ respectively in $\Ga(\Cal
T_\Phb(0,-1))$. By Proposition \ref{4.8a}, $\Bbb I_a^c\tilde\nabla_c
Y_A$ decomposes as $\tfrac{1}{2}(\nabla^HY_\Ph,\nabla^HY_\Phb)$ and
\begin{align*}
\ell^c\tilde\nabla_c Y_A=&\tfrac{1}{4}(\nabla_0
Y_\Ph-\tfrac{i}{n+2}\Rho Y_\Ph,\nabla_0 Y_\Phb+\tfrac{i}{n+2}\Rho
Y_\Phb)\\
=&\tfrac{1}{2}(iT_\al W_\Ph{}^\al+\tfrac{i}{2} SZ_\Ph,-iT_\alb
W_\Phb{}^\alb-\tfrac{i}{2} SZ_\Phb). 
\end{align*}
As observed in the proof of Proposition \ref{4.8a}, $\tilde\nabla_\k$
coincides with the double $D$--operator on weighted tractor bundles,
so by Corollary \ref{3.5} we get $\k^c\tilde\nabla_c
Y_A=\tfrac{1}{2}(-iY_\Ph,iY_\Phb)$. 

On the other hand, we have 
\begin{align*}
\Bbb Z^A{}_b=&\tW^A{}_b+\K^A\ell_b+\L^A\k_b\\
=&(\Bbb I_b^\al W^\Ph{}_\al,\Bbb I_b^\alb W^\Phb{}_\alb)+
\ell_b(iZ^\Ph,-iZ^\Phb)+\k_b\tfrac{1}{2}(iY^\Ph,-iY^\Phb).  
\end{align*}
From this, the result follows by direct evaluation. 
\end{proof}

Using this we can now directly analyse the formula \eqref{tractor-D}
for the tractor--$D$ operator.

\begin{thm*}
  Let $\tcw_{\Bbb C}[w+w']$ be a weighted complex conformal tractor
  bundle and consider the conformal tractor--$D$ operator $D_A$ which
  maps sections of $\tcw_{\Bbb C}[w+w']$ to sections of $\ce_A\otimes
  \tcw_{\Bbb C}[w+w'-1]$. Then for any $t\in\Ga(\Cal
  W(w,w'))\subset\Ga(\tcw_{\Bbb C}[w+w'])$, the holomorphic and
  anti--holomorphic parts of $D_At$ descend to sections $2D_\Ph
  t\in\Cal W(w-1,w')$ and $2D_\Phb t\in\Cal W(w,w'-1)$, and the
  operators induced in that way coincide with the CR tractor--$D$
  operators from \cite{Gover-Graham}.
\end{thm*}
\begin{proof}
  In the formula \eqref{tractor-D} from \ref{4.4}, we have to replace
  $w$ by $w+w'$, and then expand in a preferred scale. We only
  consider the holomorphic part, the anti--holomorphic part is dealt
  with in the same way. The holomorphic part of
  $(2n+2w+2w')(w+w')Y_At$ simply is $(n+w+w')(w+w')Y_\Ph t$. Next, we
  have to consider $2(n+w+w')\Bbb Z_A{}^a\tilde\nabla_a t$. Inserting
  \eqref{tWdef}, this can be written as
$$
2(n+w+w')(\tW_A{}^a+\K_A\ell^a+\L_A\k^a)\tilde\nabla_a t.
$$
Using Propositions \ref{3.6} and \ref{4.10} we conclude that the
holomorphic part of this is given by 
$$
(n+w+w')(2W_\Ph{}^\al\nabla_\al
t-iZ_\Ph(\nabla_0t+\tfrac{i(w-w')}{n+2}\Rho t)+Y_\Ph(w-w')t).
$$ 
In view of the Proposition above, the term $-(w+w')X_A\tilde\Rho t$
contributes $-(w+w')Z_\Ph\Rho t$ to the holomorphic part, and it
remains to analyse the contribution of
$-X_A\tilde\nabla^a\tilde\nabla_a t$. Using \eqref{Idef} we get 
\begin{equation}
  \label{lap1}
\tilde\nabla^a\tilde\nabla_a
t=\tilde\nabla^a(\ell_a\k^c+\k_a\ell^c+\Bbb I_a^c)\tilde\nabla_c t.   
\end{equation}
Now $\tilde\nabla^a\k_a=0$ by Proposition \ref{4.4}. The Bianchi
identity implies that
$\tilde\nabla^a\tilde\Rho_{ab}=\tilde\nabla_b\tilde\Rho$, where
$\tilde\Rho=\tilde\Rho_a{}^a$. Using this and Proposition \ref{4.4}
again, we get
$$
\tilde\nabla^a\ell_a=\tilde\nabla^a\tilde\Rho_{ab}\k^b=
\k^b\tilde\nabla_b\tilde\Rho+\tilde\Rho_{ab}\tilde\nabla^a\k^b.
$$
The first summand of this vanishes, since $\k$ is a Killing field in a
preferred scale, and the second one vanishes by symmetry of
$\tilde\Rho_{ab}$ and skew symmetry of $\tilde\nabla^a\k^b$. Hence we
can write \eqref{lap1} as
$$
\ell^a\tilde\nabla_a\k^c\tilde\nabla_c
t+\k^a\tilde\nabla_a\ell^c\tilde\nabla_c t+\tilde\nabla^a\Bbb
I_a^c\tilde\nabla_c t.
$$ 
Since the Lie bracket of $\k$ and $\ell$ vanishes, and $\k$ hooks
trivially into the conformal tractor curvature, the first two summands
are equal. Multiplying by $-X_A$ they together contribute 
$$
-i(w-w')Z_\Ph(\nabla_0t+\tfrac{i(w-w')}{n+2}\Rho t)
$$
to the holomorphic part. To analyse the last remaining term, we again
use \eqref{Idef} to get 
$$
\tilde\nabla^a\Bbb I_a^c\tilde\nabla_c
t=(\k^a\ell_b+\ell^a\k_b+\Bbb I^a_b)\tilde\nabla^b\Bbb
I_a^c\tilde\nabla_c t.
$$
Since $\k^a\Bbb I_a^c=0$ and $\ell_b\tilde\nabla^b\k^a=0$, the first
summand does not give any contribution and likewise the second summand
vanishes. Thus we are left with $\Bbb I^a_b\tilde\nabla^b\Bbb
I_a^c\tilde\nabla_c t$. But this is exactly the trace over
$\tilde\nabla^{\tilde H}\tilde\nabla^{\tilde H} t$, which descends to
the trace of $\nabla^H\nabla^H t$. This trace in turn can be written
as $\nabla^\al\nabla_\al t+\nabla^\alb\nabla_\alb t$. From  
\cite[Proposition 2.2]{Gover-Graham}, we obtain 
$$
\nabla^\alb\nabla_\alb t=\nabla^\al\nabla_\al t+\tfrac{(w-w')2(n+1)}{n+2}\Rho
t-in\nabla_0t
$$
in the case that $t$ is a density, and this formula extends to
tractors by \cite[formula (3.4)]{Gover-Graham}. 

Collecting our results, we see that the holomorphic part of $D_At$ is
given by 
\begin{align*}
  &2(n+w+w')wY_\Ph t+2(n+w+w')W_\Ph{}^\al\nabla_\al t\\
  -&2Z_\Ph(iw\nabla_0t+\nabla^\al\nabla_\al
  t+w(1+\tfrac{w'-w}{n+2})\Rho t),
\end{align*}
which is exactly twice the expression for the holomorphic CR
tractor--$D$ from \cite[section 3]{Gover-Graham}.
\end{proof}

\subsection{Almost CR-Einstein structures}\label{5.1}
Following \cite{Goalmost}, we define an \textit{almost Einstein}
structure on a pseudo-Riemannian or conformal manifold to be a
parallel section of the conformal standard tractor bundle. Almost
Einstein structures generalise the notion of Einstein manifolds since
an almost Einstein structure determines, on an open dense subset, a
scale that makes the manifold Einstein on that subset. In general, a
conformal structure does not admit an almost Einstein structure. 

Now there is an obvious analog of this condition in CR geometry.
Namely, we define an \textit{almost CR--Einstein} structure on a CR
manifold $M$ to be a parallel section of the CR standard tractor
bundle. By Proposition \ref{3.2}, there is a bijection between the
spaces of parallel CR standard tractors on $M$ and parallel conformal
standard tractors on the Fefferman space $\tilde M$. In particular,
$M$ admits an almost CR--Einstein structure if and only if $\tilde M$
admits an almost Einstein structure. 

If $I_\Phi$ is a (non-zero) parallel tractor on the CR manifold $M$
then $\sigma:=Z^\Phi I_\Phi$ is non-vanishing on an open dense
subspace of $M$; this follows because from the formula for the
connection we see that, at each point $x$, the vanishing of the
tractor covariant derivative of $I$ implies that $I$ depends only on
the 2-jet of $\sigma$ at $x$.  Away from the points where it vanishes,
$\si$ determines a scale $\si\ol{\si}$. We say that $I_\Ph$ is a
\textit{CR--Einstein} structure on $M$ if $\sigma$ is nowhere
vanishing and in this setting we will often term $\si$ itself to be a
CR-Einstein scale.

In \cite{LeePE}, J.~Lee introduced the notion of being
\textit{pseudo--Einstein} for pseudo--Hermitian structures on a CR
manifold $M$. This condition says that the Webster--Ricci tensor has
vanishing tracefree part. Motivated by ideas from tractor calculus,
F.~Leitner introduced in \cite{Leitner:TSPE} the name TSPE
(``transversally symmetric pseudo--Einstein'') for pseudo--Hermitian
structures which are pseudo--Einstein and define a transverse
symmetry. The last condition means that the Reeb vector field
corresponding to the pseudo--Hermitian structure is an infinitesimal
automorphism of the CR structure. Leitner also showed that there are
interesting relations between TSPE structures and K\"ahler--Einstein
metrics, and in particular that there are many examples of
TSPE--structures.

\begin{prop*}
  A contact form on a CR manifold $M$ induces a TSPE structure if and
  only if the corresponding scale is a CR--Einstein scale.
\end{prop*}
\begin{proof}
  Using some scale, suppose that $I_\Phi=\si
  Y_\Ph+W_\Ph{}^\al\tau_\alpha+\rho Z_\Ph$ is parallel.  From the
  formula for the connection we easily deduce that we have the
  equations
\begin{equation}\label{holoS}
\nabla_{\bar\beta} \si=0\quad \mbox{ and } \quad 
\nabla_\alpha\nabla_\beta \sigma +i\si A_{\alpha\beta}=0,
\end{equation}
which are valid in any scale and hence CR invariant. Note that, using
the formula for the covariant commutator for $f\in\ce(w,w')$
\begin{equation}\label{xcomm}
\nabla_\alpha \nabla_{\bar\beta} f- \nabla_{\bar\beta}\nabla_\alpha f
= (w-w')P_{\alpha \bar\beta} f+\frac{w-w'}{n+2}P h_{\alpha \bar\beta} f -ih_{\alpha\bar{\beta}}\nabla_0 f~,
\end{equation}
from (2.4) of \cite{Gover-Graham}, the system \nn{holoS} implies
$$
\nabla_{\bar\beta}\nabla_\alpha \si+P_{\alpha\bar\beta}
\si+h_{\alpha\bar\beta} \rho=0
$$
for a density $\rho$ in $\ce(0,-1)$.  This is another equation from
the system expressing that $I$ is covariantly parallel. 

If we now suppose that $\si$ is non-vanishing and calculate using the
pseudo-Hermitian connection $\nabla$ determined by the scale
$\si\overline{\si}$, then $\nabla_\al\si\ol{\si}=0$ and this together
with $\nabla_\al\ol{\si}=\overline{\nabla_\alb\si}$ and the above
shows that $\nabla_\beta \si =0$.  Thus the system \nn{holoS} implies
that, in the scale $\si\overline \si$, we have
$$
A_{\alpha\beta}=0 \quad \mbox{ and } \quad P_{\alpha\bar\beta}
+h_{\alpha\bar\beta} \rho/\si=0 .
$$
The first condition is well known to be equivalent to a transverse
symmetry, i.e.~the fact that the Reeb field is an infinitesimal CR
automorphisms, while the second is the equation for a pseudo-Einstein
structure in the sense of \cite{LeePE}.

Conversely suppose that we have a contact form $\th$ such that
$\Rho_{\al\beb}$ is a multiple of $h_{\al\beb}$. Now recall from
\ref{2.3} that the relation between $\th$ and the possible choices of
$\si$ such that $\si\ol{\si}$ gives the CR scale corresponding to
$\th$ is given by the fact that $\si^{n+2}\in\Ga(\ce(-n-2,0))$ is
volume normalised with respect to $\th$, and this determines
$\si^{n+2}$ up to a phase factor. Now $\ce(-n-2,0)$ is the canonical
bundle, so $\si^{n+2}$ is an $(n+1,0)$--form, and in \cite[Theorem
4.2]{LeePE} it is shown that this phase factor can be adjusted in such
a way that the resulting form is closed. This in particular implies
that $\nabla_\alb \si=0$.

Calculating in the pseudo-hermitian scale $\si\ol{\si}$ we obtain
$\nabla_\al\si=0$ as above, so if we in addition require that
$A_{\alpha\beta}=0$, we see that $\si$ solves \nn{holoS}. From
\nn{xcomm} we further get $i \nabla_0 \si =\frac{2(n+1)}{n(n+2)} P
\si$, and hence $i \nabla_\alpha \nabla_0 \si
=\frac{2(n+1)}{n(n+2)}\si \nabla_\alpha P$.

On the other hand from (2.4) of \cite{Gover-Graham} the
vanishing of $A_{\al\be}$ implies that $\nabla_\al$ and $\nabla_0$
commute on densities, so $\nabla_\alpha \nabla_{0} \si
=\nabla_{0}\nabla_\alpha \si=0$.  Thus $\nabla_\alpha P=0$. By a
similar argument we get $\nabla_{\bar\beta}\Rho =0$.  Since $P$ is a
function on $M$ it follows that $P$ is constant. Using this and the
formula for the tractor connection it is easily verified that
$$
I_\Phi:=\tfrac{1}{n+1}D_\Phi \si =\si Y_\Ph-\tfrac{1}{n}\Rho\si
Z_\Ph
$$
is annihilated by $\nabla_\alpha$ and hence by $\nabla_{\bar \beta}$
and $\nabla_0$. Thus $I_\Ph$ is parallel and since $\si$ is
non-vanishing this is a CR-Einstein structure on $M$.
\end{proof}

For the reverse implication in the proof we could also use that
\cite{Leitner:TSPE} establishes that the TSPE system implies the
Fefferman space is almost Einstein.

We have noted above that a CR--Einstein structure on $M$ induces an
almost Einstein structure on the Fefferman space $\tilde M$.
However, there is never a global Einstein metric on $\tilde M$:
\begin{thm*}
On a Fefferman space there is no Einstein metric in the conformal class.
\end{thm*}
\begin{proof}
  Suppose that $\tilde M$ is Einstein. Then it has a parallel tractor
  $I_A$ with $X^A I_A$ non--vanishing. That is writing $\si$ for the
  section $Z^\Phi I_\Phi$ of $\ce(1,0)$, we have that $\si+\overline
  \si$ is non--vanishing. By Corollary \ref{3.5} we have
  $\tilde\nabla_\k\si=i\si$ in a preferred scale. Hence under the
  fibrewise action $\rho$ of $S^1$ on $\tilde M$ we have $\rho^*_s \si
  = e^{is}\si$ and so, at any fixed point of $\tilde M$, obviously
  there is $s\in (0,2\pi] $ so that the real part of $\rho^*_s \si$
  vanishes, and this is a contradiction.
\end{proof}

\subsection*{Remarks}  (1)  The Theorem is different and
essentially stronger than that of Lee \cite{LeeF} Theorem 6.6.  In
terms of our current language the theorem of Lee shows that the
metrics determined by CR scales are never Einstein, and in this case,
there is a simpler proof: For a parallel section $I_A$ of the
complexified standard tractor bundle, the holomorphic and
anti--holomorphic parts are parallel, too. For an Einstein scale
$\al\in\Ga(\tce[1])$, the tractor $I_A= D_A\alpha$ is parallel, see
\cite{GoNur,Gau}. But for a CR scale, we have
$\al\in\ce(1/2,1/2)\subset\tce_{\Bbb C}[1]$, so the components of
$D_A\al$ are $D_\Ph\al\in\Ga(\ce(-1/2,1/2))$ and
$D_\Phb\al\in\Ga(\ce(1/2,-1/2))$ by Theorem \ref{4.10}. By Corollary
\ref{3.5}, none of the two sections can be annihilated by $\tilde
D^\nabla_{\Bbb J}$, and the latter operator coincides with
$\tilde\nabla_\k$ in a CR scale. Hence $I_A$ cannot be parallel.

\medskip

\noindent
(2) In the special case that $\dim(M)=3$, it was proved in
\cite{Lewandowski2} that there are (locally) no non-flat Einstein
metrics in the Fefferman conformal class. This follows from the
Proposition above and \cite{Leitner:TSPE}, as in this dimension the
corresponding K\"ahler-Einstein manifold has dimension 2 (and is thus flat).

\medskip

\noindent
(3) Given a CR--Einstein structure $I_\Ph$ on $M$, one can now imitate
the developments of \cite{Gopowers} to define operators with
principal part a power of the sub-Laplacian, compare these operators to
the ones obtained in \cite{Gover-Graham} and use this to prove
factorisation results. This will be taken up elsewhere.

\end{document}